\documentclass{amsart}
\usepackage{graphicx}
\usepackage{amssymb}

\numberwithin{equation}{section}

\newtheorem{thm}{Theorem}[section]
\newtheorem{lem}[thm]{Lemma}

\newtheorem{cor}[thm]{Corollary}
\newtheorem{conj}[thm]{Conjecture}
\newtheorem{defn}[thm]{Definition}
\newtheorem{obs}[thm]{Observation}
\newtheorem{prp}[thm]{Proposition}
\DeclareMathOperator{\dd}{\Delta}
\DeclareMathOperator{\lom}{\Omega^{(\ell)}}
\DeclareMathOperator{\rom}{\Omega^{(\textit{r})}}
\DeclareMathOperator{\eom}{\overline{\Omega}}
\DeclareMathOperator{\Des}{Des}
\DeclareMathOperator{\des}{des}
\DeclareMathOperator{\cDes}{cDes}
\DeclareMathOperator{\cdes}{cdes}
\DeclareMathOperator{\Pe}{Pk}
\DeclareMathOperator{\pe}{pk}
\DeclareMathOperator{\lPe}{Pk^{(\ell)}}
\DeclareMathOperator{\lpe}{pk^{(\ell)}}
\DeclareMathOperator{\rPe}{Pk^{(\textit{r})}}
\DeclareMathOperator{\rpe}{pk^{(\textit{r})}}

\DeclareMathOperator{\Sol}{Sol}
\DeclareMathOperator{\spn}{span}

\title{Enriched $P$-partitions and peak algebras}

\author{T. Kyle Petersen}

\begin{document}
\begin{abstract}
We develop a more general view of Stembridge's enriched $P$-partitions and use this theory to outline the structure of peak algebras for the symmetric group and the hyperoctahedral group. Initially we focus on commutative peak algebras, spanned by sums of permutations with the same number of peaks, where we consider several variations on the definition of ``peak." Whereas Stembridge's enriched $P$-partitions are related to quasisymmetric functions (the dual coalgebra of Solomon's type A descent algebra), our generalized enriched $P$-partitions are related to type B quasisymmetric functions (the dual coalgebra of Solomon's type B descent algebra). Using these functions, we move on to explore (non-commutative) peak algebras spanned by sums of permutations with the same set of peaks. While some of these algebras have been studied before, our approach gives explicit structure constants with a combinatorial description.
\end{abstract}

\maketitle

\section{Introduction}
This work is the result of an attempt to better understand the subalgebras of the group algebra of the symmetric group related to permutation statistics such as descent numbers and peak numbers. Much attention has been given to the so-called \emph{descent algebras}, and here we add a chapter to the story of the more recently introduced \emph{peak algebras}.

Descent algebras were first studied by Louis Solomon \cite{Solomon}, and subsequently by many others, including the papers \cite{BergeronBergeron, BergeronBergeron2, Bergeron, Cellini, Cellini2, Cellini3, Fulman, GarsiaReutenauer, Loday, Petersen}. Let $\mathfrak{S}_{n}$ be the symmetric group on $n$ elements. We think of permutations in $\mathfrak{S}_n$ as bijections \[\pi: [n] \to [n],\] where $[n]$ denotes the set $\{1,2,\ldots,n\}$. We write a permutation as an $n$-tuple, $\pi = ( \pi(1), \pi(2), \ldots, \pi(n))$. For any permutation $\pi \in \mathfrak{S}_{n}$, we say $\pi$ has a \emph{descent} in position $i$ if $\pi(i) > \pi(i+1)$. Define the set $\Des(\pi) = \{\, i \mid 1\leq i\leq n-1, \pi(i) > \pi(i+1)\,\}$ and let $\des(\pi)$ denote the number of elements in $\Des(\pi)$. We call $\Des(\pi)$ the \emph{descent set} of $\pi$, and $\des(\pi)$ the \emph{descent
number} of $\pi$. For example, the permutation $\pi = (1,4,3,2)$ has descent set
$\{2,3\}$ and descent number 2. For each subset $I$ of $[n-1]$, let \[u_{I} := \sum_{\substack{ \pi \in \mathfrak{S}_n \\ \Des(\pi) = I } } \pi \] denote the sum, in the group algebra $\mathbb{Q}[\mathfrak{S}_{n}]$, of all permutations with descent set $I$. Solomon \cite{Solomon} showed that the linear span of the $u_{I}$ forms a subalgebra of the group algebra that we call Solomon's descent algebra, denoted $\Sol(A_{n-1})$. More generally, he showed that one can define such a descent algebra, $\Sol(W)$, for any finite Coxeter group $W$.

Another sort of descent algebra, a subalgebra of $\Sol(A_{n-1})$, is given by the span of sums of permutations with the same descent number. The Eulerian descent algebra, denoted $\mathfrak{e}_n$, is defined as the linear span of the elements \[E_{i} := \sum_{ |I| = i-1} u_{I}.\] The ``Eulerian" label comes from the fact that the number of terms in $E_i$ is the Eulerian number $A_{n,i}:=\#\{\pi \in \mathfrak{S}_n \mid \des(\pi) = i-1 \}$. The Eulerian descent algebra was first studied by Jean-Louis Loday \cite{Loday} because of its use for the splitting of Hochschild homology (see also \cite{Hanlon}), and has connections with the famous card shuffling analysis of Dave Bayer and Persi Diaconis \cite{BayerDiaconis}. There are other types of Eulerian descent algebras given by considering the sums of permutations with the same number of \emph{cyclic} descents (a cyclic descent is any ordinary descent along with $n$ if $\pi(n)> \pi(1)$). While the existence of Eulerian descent algebras has not been proved for all finite Coxeter groups, thanks to Paola Cellini \cite{Cellini} the cyclic Eulerian descent algebras are known to exist in generality.\footnote{The na\"{i}ve guess at the definition of the type D Eulerian descent algebra does not work, as is easily verified for $n=3$. We make a conjecture about a general kind of Eulerian descent algebra based on Vic Reiner's \cite{Reiner} generalized $P$-partitions in section \ref{sec:conj}.} For the symmetric group and the hyperoctahedral group (to which we limit our attention in this paper), both the Eulerian and cyclic Eulerian descent algebras are well understood. See \cite{BergeronBergeron, Bergeron, Fulman}, and \cite{Petersen}. In particular, \cite{Petersen} makes use of Richard Stanley's theory of $P$-partitions (see \cite{Stanley}, chapter 4) to give a combinatorial and self-contained approach to the study of these objects. It is this thread that we pick up here to study peak algebras.

Generically, a \emph{peak} of a permutation $\pi \in \mathfrak{S}_n$ is a position $i$ such that $\pi(i-1) < \pi(i) > \pi(i+1)$, where we take $\pi(0) = \pi(n+1) = 0$. The only difference between the various types of peak sets we will study is the values of $i$ that we allow. The \emph{interior peak set} is the set of all peaks $i$ such that $1 < i < n$, the \emph{left peak set} is all peaks with $1 \leq i < n$, the \emph{right peak set} has $1 < i \leq n$, and the \emph{exterior peak set} allows $1 \leq i \leq n$. We denote these sets by $\Pe(\pi)$, $\lPe(\pi)$, $\rPe(\pi)$, and $\overline{\Pe}(\pi)$, respectively. Define the \emph{peak numbers} of $\pi$ to be the cardinalities of the peak sets: $\pe(\pi)$, $\lpe(\pi)$, $\rpe(\pi)$, and $\overline{\pe}(\pi)$. Note that $\Pe(\pi) = \lPe(\pi) \cap \rPe(\pi)$ and $\overline{\Pe}(\pi) = \lPe(\pi) \cup \rPe(\pi)$. For example, the permutation $\pi = (2,1,4,3,5)$ has $\Pe(\pi) = \{3\}$, $\lPe(\pi) = \{1,3\}$, $\rPe(\pi) = \{3,5\}$, $\overline{\Pe}(\pi) = \{1,3,5\}$, $\pe(\pi) = 1$, $\lpe(\pi) = \rpe(\pi) = 2$, and $\overline{\pe}(\pi) = 3$.

There are several relationships to be found between the algebraic structures arising from grouping permutations according to the different types of peaks. We will see that the most natural types of peaks are left peaks and interior peaks. Notice that the action of right multiplication by $\eta = (n, n-1,\ldots,1)$ on $\mathfrak{S}_n$ sends left peak numbers to right peak numbers. If $\pi = (\pi(1), \pi(2),\ldots, \pi(n))$, then $\pi\eta = (\pi(n), \pi(n-1), \ldots, \pi(1))$ and $\rPe(\pi\eta) = \{ n+1-i \,|\, i \in \lPe(\pi) \}$. In particular, $\lpe(\pi) = \rpe(\pi\eta)$. Exterior peaks are related interior peaks under left multiplication by $\eta$, with the observation that $\overline{\pe}(\pi) = \pe(\eta\pi) + 1$. These easy correspondences have some interesting consequences. We will explore connections between commutative algebras related to peaks in section \ref{sec:results}. We observe $0 \leq \pe(\pi) \leq \lfloor \frac{n-1}{2} \rfloor$, and also the number of left peaks always falls in the range $0 \leq \lpe(\pi) \leq \lfloor n/2 \rfloor$.

The study of algebras related to peaks began with John Stembridge's paper \cite{Stembridge} on \emph{enriched} $P$-partitions, followed by others, including \cite{AguiarBergeronNyman, AguiarNymanOrellana, BergeronHohlweg, BergeronMykytiukSottileWilligenburg, BilleraHsiaoWilligenburg, KrobThibon, Schocker}. While \cite{Stembridge} explores ``the algebra of peaks" related to quasisymmetric functions, it does not use enriched $P$-partitions for the study of subalgebras of $\mathbb{Q}[\mathfrak{S}_n]$ as we will here, and the only notion of peak that it uses is that of an interior peak. Kathryn Nyman \cite{Nyman} built on \cite{Stembridge} to show that there is a subalgebra of the group algebra of the symmetric group, akin to Solomon's descent algebra, formed by the linear span of \[ v_I := \sum_{\substack{\pi \in \mathfrak{S}_n \\ \Pe(\pi) = I }} \pi,\] which we call the \emph{interior peak algebra}, denoted $\mathfrak{P}_n$. Later, without the use of enriched $P$-partitions, Marcelo Aguiar, Nantel Bergeron, and Nyman \cite{AguiarBergeronNyman} showed that left peaks also give a subalgebra in this sense. We will denote the linear span of sums of permutations with the same set of left peaks by $\mathfrak{P}^{(\ell)}_n$. In \cite{AguiarBergeronNyman}, the authors also examined commutative subalgebras of the peak algebras---the ``Eulerian" peak algebras formed by sums of permutations with the same number of peaks---and showed how these subalgebras correspond to the Eulerian and cyclic Eulerian descent algebras for the hyperoctahedral group. One goal of this work is to derive some of the results of \cite{AguiarBergeronNyman} as a natural application of enriched $P$-partitions.

In this paper we will survey results of \cite{Petersen} for Eulerian descent algebras in the symmetric group and hyperoctahedral group given through the use of $P$-partitions. Subsequently we will develop a parallel study of the Eulerian peak algebras using enriched $P$-partitions. Let $\mathfrak{p}_n$, $\mathfrak{p}^{(\ell)}_n$, and $\overline{\mathfrak{p}}_n$ denote the interior, left, and exterior Eulerian peak algebras. There is no algebra generating purely by sums of permutations with the same number of right peaks.\footnote{Notice that $\eta = (n,n-1,\ldots,1)$ is the only permutation with no right peaks and $\eta^2 = (1,2,\ldots,n)$, but the identity is clearly not the only permutation with 1 right peak.} We will obtain several structure formulas that give complete sets of orthogonal idempotents for $\mathfrak{p}_n$, $\mathfrak{p}^{(\ell)}_n$, and $\overline{\mathfrak{p}}_n$, and show how these idempotents multiply with each other. In so doing, we will find that although $\mathfrak{p}_n \neq \overline{\mathfrak{p}}_n$ in general, there is a canonical isomorphism $\mathfrak{p}_n \leftrightarrow \overline{\mathfrak{p}}_n$ given by right multiplication by $\eta$. And although the linear span of sums of permutations with the same number of right peaks is not an algebra, we will see that its multiplicative closure is a commutative algebra that contains $\mathfrak{p}^{(\ell)}_n$ as a proper subalgebra. We will also show how the idempotents in the Eulerian descent algebra multiply with interior and exterior Eulerian peak idempotents and obtain a commutative peak algebra for the hyperoctahedral group.

In section \ref{sec:des} we will review some of the main ideas from \cite{Petersen}, since it is this approach that we will mirror later on. In section \ref{sec:results} we will give an overview of the results for the Eulerian peak algebras, though the theorems come well before their proofs or even a precise statement of some of the key definitions---these are provided in section \ref{sec:epp}, which gives a rigorous treatment of enriched $P$-partitions. Section \ref{sec:peak} presents the proofs of the results of section \ref{sec:results} using the theory of enriched $P$-partitions.

A topic not explored deeply in this paper is that of further applications of enriched $P$-partitions through their quasisymmetric generating functions. Ira Gessel \cite{Gessel} outlined how generating functions for ordinary $P$-partitions give a natural basis for the space of quasisymmetric functions. Stembridge \cite{Stembridge} defined generating functions for enriched $P$-partitions that form a subring of the ring of quasisymmetric functions. Just as Stembridge's enriched $P$-partitions connect with quasisymmetric functions (the dual coalgebra to $\Sol(A_{n-1})$), the new types of enriched $P$-partitions we present here connect to type B quasisymmetric functions (the dual coalgebra to $\Sol(B_n)$), as studied by Chak-On Chow \cite{Chow} using ordinary type B $P$-partitions. In section \ref{sec:qsym}, we will lay the groundwork for future work in this direction. One significant result we establish is the existence of a (noncommutative) type B peak algebra, $\mathfrak{P}_{B,n}$, spanned by sums of signed permutations $\pi$ with the same peak set and the same sign on $\pi(1)$. Through our approach we will obtain combinatorially described structure constants for each of the algebras $\mathfrak{P}_n$, $\mathfrak{P}^{(\ell)}_n$, and $\mathfrak{P}_{B,n}$.

\section{$P$-partitions and Eulerian descent algebras}\label{sec:des}

In this section we review the prototypical example of the $P$-partition approach, since we will model our work with enriched $P$-partitions on this approach. We follow \cite{Petersen}, whose presentation was directly motivated by Gessel's paper \cite{Gessel}. We will also describe the type B results and make a general conjecture for finite Coxeter groups.

\subsection{Type A}

\begin{figure} [h]
\centering
\includegraphics{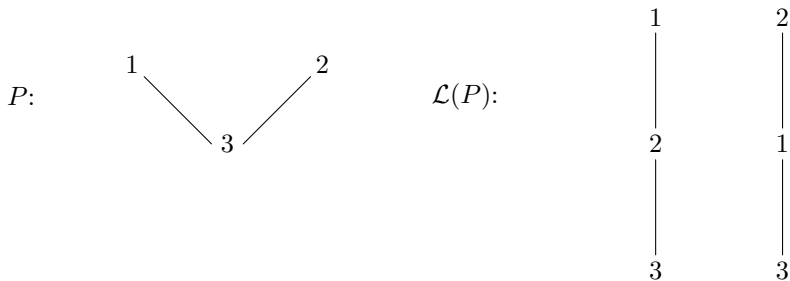}
\caption{Linear extensions of a poset $P$.\label{fig:hasse}}
\end{figure}

The ``$P$" in $P$-partition stands for a partially ordered set, or poset. For our purposes, we assume that all posets $P$, with partial order $<_{P}$, are finite. And unless otherwise noted, if $|P| = n$, then the elements of $P$ are labeled distinctly with the numbers $1,2,\ldots,n$. We will sometimes describe a poset by its Hasse diagram, as in Figure \ref{fig:hasse}.

\begin{defn}[$P$-partition]\label{def:p}
Let $X = \{x_{1},x_{2},\ldots \}$ be a countable, totally ordered
set. For a given poset $P$, a \emph{$P$-partition} is an order-preserving\footnote{We note that this definition differs from Richard Stanley's \cite{Stanley} in that our maps are order-preserving, whereas his are order reversing, i.e., $f(i) \geq f(j)$ if $i <_P j$. He calls the maps of Definition \ref{def:p} (perhaps misleadingly) \emph{reverse} $P$-partitions.} map $f: [n] \to X$ such that:
\begin{enumerate}
\item $f(i) \leq f(j)$ if $i <_{P} j$

\item $f(i) < f(j)$ if $i <_{P} j$ and $i > j$ in $\mathbb{Z}$
\end{enumerate}
\end{defn}
For our purposes we usually think of $X$ as a subset of the
positive integers. Let $\mathcal{A}(P)$ denote the set of all
$P$-partitions. When $X$ has finite cardinality $k$, then the number of
$P$-partitions must also be finite. In this case, define the \emph{order
polynomial}, denoted $\Omega(P;k)$, to be the number of
$P$-partitions $f:[n] \to X$. Order polynomials play a critical role in understanding the Eulerian descent algebra.

We can think of any permutation $\pi \in \mathfrak{S}_{n}$ as a poset with the total order $\pi(s) <_{\pi} \pi(s+1)$. For example, the permutation $\pi = (3,2,1,4)$ has $3 <_{\pi} 2
<_{\pi} 1 <_{\pi} 4$ as a poset. With this convention, the set
of all $\pi$-partitions is easily characterized in terms of descents. Observe that
$\mathcal{A}(\pi)$ is the set of all functions $f: [n] \to X$ such
that \[f(\pi(1)) \leq f(\pi(2)) \leq \cdots \leq f(\pi(n)),\]
and whenever $\pi(s) > \pi(s+1)$ (i.e., $s \in \Des(\pi)$), then $f(\pi(s))< f(\pi(s+1))$. The set of all $\pi$-partitions where $\pi = (3,2,1,4)$ is all maps $f$ such that $f(3) < f(2) < f(1) \leq f(4)$.

For a poset $P$ with $n$ elements, let $\mathcal{L}(P)$ denote its Jordan-H\"{o}lder set: the set of all
permutations of $[n]$ which extend $P$ to a total order. This set is also called the set
of ``linear extensions" of $P$. For example let $P$ be the poset
defined by $1 >_{P} 3 <_{P} 2$. In linearizing $P$ we form a
total order by retaining all the relations of $P$ but introducing
new relations so that any element is comparable to any other. In
this case, 1 and 2 are not comparable, so we have exactly two ways
of linearizing $P$: $3 < 2 < 1$ or $3 < 1 < 2$. These correspond
to the permutations $(3,2,1)$ and $(3,1,2)$. Let us make the
following observation.
\begin{obs}\label{ob1} A permutation $\pi$ is in $\mathcal{L}(P)$ if and only if $i<_{P}
j$ implies $\pi^{-1}(i) < \pi^{-1}(j)$.
\end{obs}
In other words, if $i$ is ``below" $j$ in the Hasse diagram of the poset $P$, it had better be below $j$ in any linear extension of the poset.
We now prove what we call the fundamental lemma of $P$-partitions.
\begin{lem}
The set of all $P$-partitions of a poset $P$ is the disjoint union
of the set of $\pi$-partitions of all linear extensions $\pi$ of
$P$: \[\mathcal{A}(P) = \coprod_{\pi \in \mathcal{L}(P)}
\mathcal{A}(\pi).\]
\end{lem}
\begin{proof}
The proof follows from induction on the number of incomparable
pairs of elements of $P$. If there are no incomparable pairs, then
$P$ has a total order and already represents a permutation.
Suppose $i$ and $j$ are incomparable in $P$. Let $P_{ij}$ be the
poset formed from $P$ by introducing the relation $i < j$. Then it
is clear that $\mathcal{A}(P) = \mathcal{A}(P_{ij}) \coprod
\mathcal{A}(P_{ji})$. We continue to split these posets (each with strictly fewer incomparable pairs)
until we have a collection of totally ordered chains corresponding to distinct linear
extensions of $P$.
\end{proof}
\begin{cor}The order polynomial of a poset is the sum of the order polynomials of its linear extensions:
\[\Omega(P;k) = \sum_{\pi \in \mathcal{L}(P)} \Omega(\pi;k).\]
\end{cor}

The fundamental lemma tells us that in order to study $P$-partitions, we can focus on the the case where $P$ is a totally ordered chain---a more straightforward task. In particular, counting $\pi$-partitions that map into a finite set is not too difficult. Notice that for any permutation $\pi$ and any positive integer $k$,
\begin{equation}\label{eq:binom}
\begin{aligned}
\Omega(\pi;k) & := \#\{ f: [n] \to [k] \,|\, 1 \leq f(\pi(1)) \leq f(\pi(2)) \leq \cdots \leq f(\pi(n)) \leq k \\
 & \qquad \qquad \qquad \qquad \qquad \mbox{ and } f(\pi(s)) < f(\pi(s+1)) \mbox{ if } s \in \Des(\pi) \}  \\
 & =  \#\{ (i_1, i_2, \ldots, i_n) \in \mathbb{Z}^n \,|\, 1 \leq i_1 \leq i_2 \leq \cdots \leq i_n \leq k \mbox{ and } i_s < i_{s+1} \mbox{ if } s \in \Des(\pi) \} \\
 & =  \#\{ (i_1, i_2, \ldots, i_n) \in \mathbb{Z}^n \,|\, 1 \leq i_1 < i_2 < \cdots < i_n \leq k +n -1 - \des(\pi) \}  \\
 & =  \binom{ k+ n-1 - \des(\pi)}{n}.
\end{aligned}
\end{equation}
The third equality above is given by the observation that if $i \leq j$, then $i < j+1$. For example, the number of solutions to $1 \leq i_1 < i_2 \leq i_3 \leq 4$ is the same as the number of solutions to $1 \leq i_1 < i_2 < i_3 + 1 \leq 4+1$ or the solutions to $1 \leq i_1 < i_2 < i'_3 \leq 5$. As immediate consequences of \eqref{eq:binom}, we see that the order polynomial of a permutation is a polynomial of degree $n$ with no constant term and depends only on the number of descents of the permutation. For fixed $n$, let $\Omega(i;x) = \binom{x + n - i}{n}$ denote the order polynomial for any permutation of $[n]$ with $i-1$ descents.

Define the following element of $\mathbb{Q}[\mathfrak{S}_n][x]$, which we refer to as the ``structure polynomial": \[ \phi(x) := \sum_{\pi\in\mathfrak{S}_{n}}
\Omega(\pi;x) \pi = \sum_{i=1}^{n}\Omega(i;x)E_{i},\] where we recall that $E_i$ is the sum of all permutations with $i-1$ descents. By construction, the structure polynomial is a polynomial in $x$, with coefficients in the group algebra $\mathbb{Q}[\mathfrak{S}_n]$, of degree $n$ and with no constant term. In other words, it has exactly as many nonzero terms as there are possible descent numbers. The Eulerian subalgebra $\mathfrak{e}_n$ is described by the following multiplication of structure polynomials, which is essentially an unsigned version of Loday's Th\'{e}or\`{e}me 1.7 (see \cite{Loday}).
\begin{thm}[Gessel \cite{Gessel}]\label{thm:ges}
As polynomials in $x$ and $y$ with coefficients in the group algebra, we have
\begin{equation}\label{eq:ges}
\phi(x)\phi(y) =\phi(xy).
\end{equation}
\end{thm}
Define elements $e_{i}$ in the group algebra by $\displaystyle \phi(x) = \sum_{i=1}^{n} e_{i}x^{i}$. By examining the coefficients of $x^i y^j$ in \eqref{eq:ges}, it is clear that the $e_{i}$ are orthogonal idempotents: $e_{i} e_j = \delta_{ij} e_{i}$, where $\delta_{ij}$ is Kronecker's delta function: $\delta_{ij} = 1$ if $i=j$, $\delta_{ij} = 0$ otherwise. The following result is well established.

\begin{cor}\label{cor:desalg}
The Eulerian descent algebra $\mathfrak{e}_n$ is commutative of dimension $n$.
\end{cor}

\begin{proof}[Proof of Corollary \ref{cor:desalg}]
It will suffice to show that \[\spn\{E_1, \ldots, E_n\} = \spn\{e_1, \ldots, e_n\},\] since the $e_i$ obviously span an $n$-dimensional, commutative algebra. It is immediate from the definition of the structure polynomial that $\spn\{e_i\} \subset \spn\{E_i\}$. We need only show the reverse inclusion.

First, notice that $\phi(1) = E_1 = e_1 + \cdots + e_n$, since $\Omega(\pi;1) = \binom{n - \des(\pi)}{n}$ is zero if $\des(\pi) > 0$. We now proceed by induction. Suppose that for some $k \leq n$, all $E_i$, $i = 1,\ldots, k-1$ can be written as linear combinations of the $e_i$. Then we have
\begin{align*}
\phi(k) & =  \sum_{i=1}^n \binom{n+k-i}{n} E_i \\
 & =  \binom{n+k-1}{n}E_1 + \cdots + \binom{n+1}{n}E_{k-1} + E_k \\
 & =  \sum_{i=1}^{n} k^i e_i.
\end{align*}
By the induction hypothesis, the expression \[E_k = \sum_{i=1}^{n} k^i e_i - \sum_{i=1}^{k-1} \binom{n+k-i}{n}E_i,\] is in $\spn\{e_i\}$, and the result follows.
\end{proof}

We could refer the reader to \cite{Petersen} for proof of Theorem \ref{thm:ges}, but because later arguments are so similar, we include it below as a kind of ground-level proof on which later proofs are built. We point out that in order to prove that the formulas in this paper hold as polynomials in $x$ and $y$, it will suffice to prove that they hold for all pairs of positive integers. It is not hard to verify this fact, and it is implicit in many of the proofs presented in this paper.

\begin{proof}[Proof of Theorem \ref{thm:ges}]
If we write out $\phi(xy)=\phi(x)\phi(y)$ using the
definition of the structure polynomials, we have
\begin{align*}
\sum_{\pi\in\mathfrak{S}_{n}}\Omega(\pi;xy)\pi & =
\sum_{\sigma\in\mathfrak{S}_{n}}\Omega(\sigma;x) \sigma
\sum_{\tau\in\mathfrak{S}_{n}}\Omega(\tau;y)\tau \\
 & =  \sum_{\sigma,\tau \in \mathfrak{S}_n }\Omega(\sigma;x)\Omega(\tau;y) \sigma\tau.
\end{align*}
If we equate the coefficients of $\pi$ we have
\begin{equation}\label{eq:order}
\Omega(\pi;xy) = \sum_{\sigma\tau = \pi}
\Omega(\sigma;x) \Omega(\tau;y).
\end{equation}
Clearly, if formula \eqref{eq:order} holds for all $\pi$, then formula
\eqref{eq:ges} is true. Let $x =k$ and $y=l$ be positive integers and consider the left hand
side of equation \eqref{eq:order}. To compute the order polynomial $\Omega(\pi;kl)$ we need to count the number of $\pi$-partitions $f:[n] \to X$, where $X$ is some totally ordered set with $kl$ elements. But instead of using $[kl]$ as our image set, we will use a different totally ordered set of the same cardinality. Let us count the $\pi$-partitions $f:[n]\to [l]\times[k]$. This is equal to the number of solutions to
\begin{equation}
(1,1) \leq (i_1,j_1) \leq (i_2,j_2) \leq \cdots \leq (i_n,j_n)
\leq (l,k)
\mbox{ and } (i_{s},j_s) < (i_{s+1},j_{s+1}) \mbox{ if } s \in \Des(\pi). \label{eq:lex}
\end{equation}
Here we take the \emph{lexicographic
ordering} on pairs of integers. Specifically, $(i,j) < (i',j')$ if
$i < i'$ or else if $i = i'$ and $j < j'$.

To get the result we desire, we will sort the set of all solutions to \eqref{eq:lex} into distinct cases indexed by subsets $I \subset [n-1]$. The sorting depends on $\pi$ and proceeds as follows. Let $F = ( (i_1, j_1), \ldots, (i_n, j_n) )$ be any solution to \eqref{eq:lex}. For any $s=1,2,\ldots,n-1$, if $\pi(s) < \pi(s+1)$, then $(i_{s},j_{s}) \leq (i_{s+1},j_{s+1})$, which falls into one of two mutually exclusive cases:
\begin{align}
i_{s} \leq i_{s+1} & \mbox{ and }  j_{s}\leq j_{s+1}, \mbox{ or } \label{eq:split1}\\
i_{s} < i_{s+1} & \mbox{ and }  j_{s} > j_{s+1}. \label{eq:split2}
\end{align}
If $\pi(s) > \pi(s+1)$, then $(i_{s},j_{s}) < (i_{s+1},j_{s+1})$,
which means either:
\begin{align}
i_{s} \leq i_{s+1} & \mbox{ and }  j_{s} < j_{s+1}, \mbox{ or } \label{eq:split3}\\
i_{s} < i_{s+1} & \mbox{ and }  j_{s} \geq j_{s+1},\label{eq:split4}
\end{align}
also mutually exclusive. Define $I_F = \{ s \in [n-1]\setminus \Des(\pi) \,\mid\, j_s > j_{s+1} \} \cup \{ s \in \Des(\pi) \,\mid\, j_s \geq j_{s+1} \}$. Then $I_F$ is the set of all $s$ such that either \eqref{eq:split2} or \eqref{eq:split4} holds for $F$. Notice that in both cases, $i_s < i_{s+1}$. Now for any $I \subset [n-1]$, let $S_I$ be the set of all solutions $F$ to \eqref{eq:lex} satisfying $I_F = I$. We have split the solutions of \eqref{eq:lex} into $2^{n-1}$ distinct cases indexed by all the different subsets $I$ of $[n-1]$.

Say $\pi = (2,1,3)$. Then we want to count the number of solutions to \[(1,1)\leq (i_{1},j_{1}) < (i_{2},j_{2}) \leq (i_{3},j_{3}) \leq (l,k)\] which splits into four distinct cases:

\begin{table}[h]
\begin{tabular}{|c|c c c|}
\hline & & & \\
$\emptyset$ & $i_{1} \leq i_{2} \leq i_{3}$ & and & $j_{1} < j_{2} \leq j_{3}$ \\
& & & \\
$\{\,1\,\}$ & $i_{1} < i_{2} \leq i_3$ & and & $j_1 \geq j_2 \leq j_3$ \\
& & &\\
$\{\,2\,\}$ & $i_1 \leq i_2 < i_3$ & and & $j_1 < j_2 > j_3$ \\
& & & \\
$\{\,1,2\,\}$ & $i_1 < i_2 < i_3$ & and & $j_1 \geq j_2 > j_3$\\
& & & \\
\hline
\end{tabular}
\end{table}

We now want to count all the solutions contained in each of these cases and add them up. For a fixed subset $I$ we will use the theory of $P$-partitions to count the number of solutions for the set of inequalities first for the $j_{s}$'s and then for the $i_{s}$'s. Multiplying will give us the number of solutions in $S_I$; we do the same for the remaining subsets and sum to obtain the final result. For $I = \{\,1\,\}$ in the example above, we would count first the number of integer solutions to $j_1 \geq j_2 \leq j_3$, with $1\leq j_s \leq k$, and then we multiply this number by the number of solutions to $1 \leq i_1 < i_2 \leq i_3 \leq l$ to obtain the cardinality of $S_{\{1\}}$.  We will now carry out the computation in general.

For any particular $I\subset [n-1]$, form the poset $P_{I}$ of the elements
$1,2,\ldots,n$ by $\pi(s) <_{P_{I}} \pi(s+1)$ if $s \notin I$,
$\pi(s) >_{P_{I}} \pi(s+1)$ if $s\in I$. We form a ``zig-zag" poset of $n$ elements labeled consecutively by $\pi(1),
\pi(2),\ldots,\pi(n)$ with downward zigs corresponding to the
elements of $I$. For example, if $I =\{2,3\}$ for $n=5$, then $P_{I}$ has $\pi(1) <
\pi(2) > \pi(3) > \pi(4) < \pi(5)$. See Figure \ref{fig:zigzag}.

\begin{figure} [h]
\centering
\includegraphics{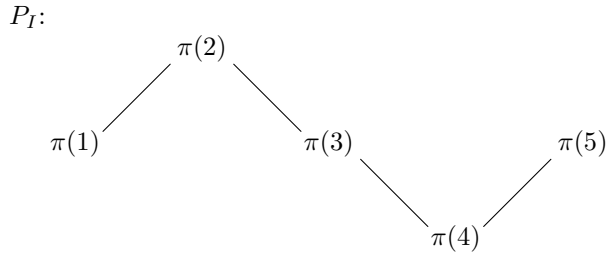}
\caption{The ``zig-zag" poset $P_{I}$ for $I = \{2,3\} \subset
[5]$. \label{fig:zigzag}}
\end{figure}

For any solution in $S_I$, let $f: [n] \to [k]$ be defined by $f(\pi(s)) = j_{s}$ for $1\leq
s\leq n$. We will show that $f$ is a $P_{I}$-partition. If $\pi(s) <_{P_{I}} \pi(s+1)$ and $\pi(s) < \pi(s+1)$
in $\mathbb{Z}$, then \eqref{eq:split1} tells us that $f(\pi(s)) =
j_{s} \leq j_{s+1} = f(\pi(s+1))$. If $\pi(s) <_{P_{I}} \pi(s+1)$
and $\pi(s) > \pi(s+1)$ in $\mathbb{Z}$, then \eqref{eq:split3} tells
us that $f(\pi(s)) = j_{s} < j_{s+1} = f(\pi(s+1))$. If $\pi(s)
>_{P_{I}} \pi(s+1)$ and $\pi(s) < \pi(s+1)$ in $\mathbb{Z}$,
then \eqref{eq:split2} gives us that $f(\pi(s))= j_{s} > j_{s+1} =
f(\pi(s+1))$. If $\pi(s) >_{P_{I}} \pi(s+1)$ and $\pi(s) >
\pi(s+1)$ in $\mathbb{Z}$, then \eqref{eq:split4} gives us that
$f(\pi(s)) = j_{s} \geq j_{s+1} = f(\pi(s+1))$. In other words, we
have verified that $f$ is a $P_{I}$-partition. So for any particular solution in $S_I$, the $n$-tuple $(j_1, \ldots, j_n)$ can be thought of as a $P_{I}$-partition. Conversely, any $P_{I}$-partition $f$ gives a solution in $S_I$ since if $j_s = f(\pi(s))$, then $(( i_1,j_1),\ldots,(i_n,j_n)) \in S_I$ if and only if $1 \leq i_1 \leq \cdots \leq i_n \leq l$ and $i_s < i_{s+1}$ for all $s \in I$. We can therefore turn our attention to counting $P_{I}$-partitions.

Let $\sigma \in \mathcal{L}(P_{I})$. Then for any
$\sigma$-partition $f$, we get \[1\leq f(\sigma(1)) \leq f(\sigma(2))\leq \cdots \leq f(\sigma(n)) \leq
k\] with $f(\sigma(s)) < f(\sigma(s+1))$ if $s \in \Des(\sigma)$.
The number of solutions to this set of inequalities is by definition $\Omega(\sigma;k)$.

Recall by Observation \ref{ob1} that $\sigma^{-1}\pi(s) < \sigma^{-1}\pi(s+1)$ if $\pi(s) <_{P_{I}} \pi(s+1)$, i.e., if $s \notin I$. If $\pi(s) >_{P_{I}} \pi(s+1)$ then $\sigma^{-1}\pi(s) > \sigma^{-1}\pi(s+1)$ and $s \in I$. We get that $\Des(\sigma^{-1}\pi) = I$ if and only if $\sigma \in \mathcal{L}(P_I)$. Set $\tau = \sigma^{-1}\pi$. The number of solutions to \[1 \leq i_{1} \leq \cdots \leq i_{n} \leq l \mbox{ and } i_{s} < i_{s+1} \mbox{ if } s\in \Des(\tau)\] is defined to be $\Omega(\tau;l)$. Now for a given $I$, the number of solutions in $S_I$ is \[\sum_{\substack{\sigma\in\mathcal{L}(P_{I}) \\ \sigma\tau = \pi}} \Omega(\sigma;k)\Omega(\tau;l). \] Summing over all subsets $I \subset [n-1]$, we
can write the number of all solutions to \eqref{eq:lex} as \[\sum_{\sigma\tau=\pi}\Omega(\sigma;k)\Omega(\tau;l),\] and so we have derived formula \eqref{eq:order}.
\end{proof}

We will now present the main results from \cite{Petersen}, though Theorem \ref{thm:chow} first appeared in \cite{Chow}. We omit the proofs here, but remark that conceptually they follow the same lines as the proof of Theorem \ref{thm:ges} above.

For any $\pi \in \mathfrak{S}_n$, define a \emph{cyclic descent} to be any $i=1,2,\ldots,n-1$ such that $\pi(i) > \pi(i+1)$ along with $i=n$ if $\pi(n) > \pi(1)$. Define $\cDes(\pi)$ to be the \emph{cyclic descent set} of $\pi$, and $\cdes(\pi)$ to be the \emph{cyclic descent number}. Let $E_{i}^{(c)}$ be the sum in the group algebra of all those permutations with $i$ cyclic descents, and let $\mathfrak{e}^{(c)}_n$ denote the linear span of the $E_i^{(c)}$, called the \emph{cyclic} Eulerian descent algebra. Define $\Omega^{(c)}(\pi;x) = \frac{1}{n}\binom{x+n-1-\cdes(\pi)}{n-1}$. We write the cyclic structure polynomial \[\varphi(x) := \sum_{\pi\in\mathfrak{S}_{n}}\Omega^{(c)}(\pi;x)\pi =
\sum_{i=1}^{n-1}\Omega^{(c)}(i;x) E_{i}^{(c)} = \sum_{i=1}^{n-1} e_i^{(c)} x^i,\] and the structure of $\mathfrak{e}^{(c)}_n$ is given by the following theorem.
\begin{thm}[Petersen \cite{Petersen}]\label{thmcyc}
As polynomials in $x$ and $y$ with coefficients in the group
algebra of the symmetric group, we have
\[\varphi(x)\varphi(y) = \varphi(xy).\]
\end{thm}

Thus the elements $e_i^{(c)}$ are orthogonal idempotents.

\begin{cor}\label{cor:cdesalg}
The cyclic Eulerian descent algebra $\mathfrak{e}^{(c)}_n$ is commutative of dimension $n-1$.
\end{cor}

Further, it is not too much work to see the Eulerian descent algebra $\mathfrak{e}_{n-1}$ is isomorphic to the cyclic Eulerian descent algebra $\mathfrak{e}^{(c)}_n$ via the map \[ \pi \mapsto \sum_{i = 1}^n \widehat{\pi}\omega^i,\] where for any $\pi \in \mathfrak{S}_{n-1}$, $\widehat{\pi} = (\pi(1), \pi(2), \ldots, \pi(n-1), n)$, and $\omega = (2, 3,\ldots, n, 1)$ is the $n$-cycle.

\subsection{Type B}

For the hyperoctahedral group $\mathfrak{B}_n$, the group of signed permutations, the results are similar. Let $\pm[n]$ denote the set $\{-n, -n+1, \ldots, -1,0,1,\ldots, n-1, n\}$. We think of signed permutations as bijections $\pi: \pm[n] \to \pm[n]$ with the property that $\pi(-i) = -\pi(i)$. Thus, $\pi(0) = 0$, and the image of $\pi$ is determined by $\pi(1), \pi(2), \ldots, \pi(n)$. When considering descents in the hyperoctahedral group, the only difference from the symmetric group is that we need to allow a descent at the beginning of the permutation. We define the \emph{descent set} $\Des_B(\pi)$ of a signed permutation to be the set of all $i \in [0,n-1] = \{0,1,2,\ldots,n-1\}$ such that $\pi(i) > \pi(i+1)$. So in particular if $\pi(1)$ is negative, then $0$ is a descent of $\pi$. The \emph{descent number} of $\pi$ is denoted $\des_B(\pi)$ and is equal to the cardinality of $\Des_B(\pi)$. The \emph{cyclic descent set}, $\cDes_B(\pi)$, is the set of all ordinary descents along with $n$ if $\pi(n) > \pi(0) = 0$, and we denote the \emph{cyclic descent number} by $\cdes_B(\pi)$. As a simple example, the signed permutation $(-2,1)$ has descent set $\{0\}$, cyclic descent set $\{0,2\}$, descent number $1$, and cyclic descent number $2$.

For fixed $n$, we let $E_{B,i}$ be the sum of all signed permutations with $i-1$ descents, and $E^{(c)}_{B,i}$ be the sum of all signed permutations with $i$ cyclic descents. The Eulerian and cyclic Eulerian descent algebras are denoted $\mathfrak{e}_{B,n}$ and $\mathfrak{e}^{(c)}_{B,n}$. For the hyperoctahedral group, the order polynomials turn out to be $\Omega_B (\pi;x) = \binom{x+n-\des_B(\pi)}{n}$ in the ordinary case, and $\Omega^{(c)}_B (\pi; x) = \binom{x+n-\cdes_B(\pi)}{n}$ in the cyclic case. If we define \[\phi_B(x) := \sum_{\pi\in\mathfrak{B}_{n}}\Omega_B(\pi;(x-1)/2) \pi =
\sum_{i=1}^{n+1} \Omega_B(i;(x-1)/2) E_{B,i} = \sum_{i=0}^{n} e_i x^i, \] (note the constant term) and \[ \varphi_B(x) := \sum_{\pi\in\mathfrak{B}_{n}} \Omega^{(c)}_B(\pi;x/2) \pi = \sum_{i=1}^{n} \Omega^{(c)}_B(i;x/2) E^{(c)}_{B,i} = \sum_{i=1}^n e_i^{(c)} x^i, \] we have the following theorems.

\begin{thm}[Chow \cite{Chow}]\label{thm:chow}
As polynomials in $x$ and $y$ with coefficients in the group algebra of the hyperoctahedral group, we have
\[\phi_B(x)\phi_B(y) = \phi_B(xy). \]
\end{thm}

\begin{cor}\label{cor:desBalg}
The type B Eulerian descent algebra $\mathfrak{e}_{B,n}$ is commutative of dimension $n+1$.
\end{cor}

\begin{thm}[Petersen \cite{Petersen}]\label{thm:cyclicB}
As polynomials in $x$ and $y$ with coefficients in the group
algebra of the hyperoctahedral group we have \[\varphi_B(x)\varphi_B(y) = \varphi_B(xy).\]
\end{thm}

\begin{cor}\label{cor:cdesBalg}
The Eulerian descent algebra $\mathfrak{e}_n$ is commutative of dimension $n$.
\end{cor}

\begin{thm}[Petersen \cite{Petersen}]\label{thm:idealB}
As polynomials in $x$ and $y$ with coefficients in the group algebra of the hyperoctahedral group we have \[\varphi_B(x)\phi_B(y) = \phi_B(y)\varphi_B(x) = \varphi_B(xy).\]
\end{thm}

Thus we have the following multiplication rule: $e_i^{(c)} e_j = e_j e_i^{(c)} =  \delta_{ij} e_i^{(c)}$. Let $\dot{\mathfrak{e}}_{B,n}$ denote the algebra formed by the span of both the Eulerian elements, $E_{B,i}$, and cyclic Eulerian elements, $E_{B,i}^{(c)}$.

\begin{cor}\label{cor:idealB}
The algebra $\dot{\mathfrak{e}}_{B,n}$ is commutative of dimension $2n$.
\end{cor}

\begin{proof}
Commutativity is obvious; it may not be obvious that the dimension is $2n$ (rather than $2n+1$).
For $i = 1,2,\ldots,n$, let $F_{i}^{-}$ be the sum of all signed permutations with $i$ cyclic descents and $\pi(n) <0$, let $F_{i}^{+}$ be the sum of all signed permutations with $i$ cyclic descents and $\pi(n) >0$. Then
\begin{enumerate}
\item $E_{B,1} =  F_{1}^{+}$,
\item $E_{B,n+1} =  F_{n}^{-}$,
\item $E_{B,i} = F_{i-1}^{-} + F_{i}^{+}$ for  $1 < i < n+1$, and
\item $E_{B,i}^{(c)} = F_{i}^{-} + F_{i}^{+}$ for $1 \leq i \leq n$.
\end{enumerate}
Then we see that the $F_{i}^{+}$, $F_{i}^{-}$, which are obviously linearly independent, span the $E_{B,i}$, $E_{B,i}^{(c)}$.
\end{proof}

Theorems \ref{thm:chow} and \ref{thm:cyclicB} tell us that $e_i$, $e^{(c)}_i$, are orthogonal idempotents for the type B Eulerian and cyclic Eulerian descents algebras respectively. But Theorem \ref{thm:idealB} presents something not seen in the case of the symmetric group algebra. It allows us to conclude that in $\dot{\mathfrak{e}}_{B,n}$, the cyclic Eulerian subalgebra $\mathfrak{e}^{(c)}_{B,n}$ is an ideal. This fact was observed in \cite{AguiarBergeronNyman} and a similar result will be seen for the peak algebras of the symmetric group algebra.

\subsection{Conjecture: a general commutative descent algebra}\label{sec:conj}

As mentioned in the introduction, Cellini \cite{Cellini} established that the cyclic Eulerian descent algebras exist for all finite Coxeter groups $W$, whereas an ordinary Eulerian descent algebra (spanned by sums of permutations with the same number of descents) fails even for type D. However, we conjecture that order polynomials give a way to construct a general commutative subalgebra of $\mathbb{Q}[W]$, for any finite Coxeter group $W$. Phrased in terms of root systems, Reiner \cite{Reiner2} has defined $P$-partitions and order polynomials for any $W$. His definition coincides with the definitions for type A and type B. We make the following conjecture.

\begin{conj}\label{conj:des}
Let $W$ be any finite Coxeter group. For any $w \in W$, we have
\[ \Omega_W(w;t_W(xy)) = \sum_{uv = w} \Omega_W(u;t_W(x)) \Omega_W(v;t_W(y)),\]
where $t_W$ is a linear function that depends only on $W$.
\end{conj}

Indeed, this formula is true for $A_n$, $B_n$, and (according to Chow \cite{Chow2}) for $D_n$, where $t_{A}(x) = x$, $t_{B}(x) = t_{D}(x) = (x-1)/2$. If Conjecture \ref{conj:des} is true, then we can define the following polynomial in $\mathbb{Q}[W][x]$,
\[ \phi_W(x) = \sum_{w \in W} \Omega_W (w;t_W(x)) w, \]
which would satisfy
\[ \phi_W(x)\phi_W(y) = \phi_W(xy).\] Then the coefficients of $\phi_W(x)$ would give orthogonal idempotents for a commutative subalgebra of $\mathbb{Q}[W]$. The order polynomials in types A and B depend only on the number of descents. In type D, \cite{Chow2} shows that the order polynomial depends on the number of descents of $w$ and $\overline{w}$, where $\overline{w}(1) = - w(1)$, and $\overline{w}(i) = w(i)$ for $i \geq 2$.

The remaining cases left to prove for Conjecture \ref{conj:des} should be straightforward to verify, if a little tedious. Is there a case-independent proof?

\section{Commutative peak algebras}\label{sec:results}

The results of this section establish the structure and interactions of the Eulerian peak algebras in the symmetric group algebra, and shows how the interior and exterior peak algebras relate to Eulerian descent algebras. Further, we obtain the idempotents for a type B Eulerian peak algebra. The main tools we use come from the theory of \emph{enriched} $P$-partitions, which we develop rigorously in section \ref{sec:typeAepp}. For now let us take some facts for granted.

\subsection{Type A}

To any poset $P$ we will associate a polynomial $\Omega'(P;x)$, that we call the \emph{enriched order polynomial}. In section \ref{sec:typeAepp} we will prove that for a permutation $\pi$, $\Omega'(\pi;x/2)$ is an even or odd polynomial of degree $n$, with no constant term, that depends only on the number of interior peaks of $\pi$. Similarly, we have what is called the \emph{left enriched order polynomial}, denoted $\lom(P;x)$. The pertinent fact is that $\lom(\pi;(x-1)/2)$ is an even or odd polynomial of degree $n$ that depends only on the number of left peaks of $\pi$.  We form the structure polynomials
\begin{align*}
 \rho(x) & :=  \sum_{\pi \in \mathfrak{S}_n} \Omega'(\pi;x/2) \pi =  \sum_{i=1}^{\lfloor \frac{n+1}{2} \rfloor} \Omega'(i;x/2)E'_{i},\\
 & =
\begin{cases}
\displaystyle \sum_{i=1}^{n/2} e'_i x^{2i} & \mbox{ if $n$ is even, }\\
\displaystyle \sum_{i=1}^{(n+1)/2} e'_i x^{2i-1} & \mbox{ if $n$ is odd, }
\end{cases}
\end{align*}
and
\begin{align*}
\rho^{(\ell)}(x) & := \sum_{\pi \in \mathfrak{S}_n} \lom(\pi;(x-1)/2) \pi =  \sum_{i=1}^{\lfloor \frac{n}{2} \rfloor +1} \lom(i;(x-1)/2)E^{(\ell)}_{i},\\
 & =
\begin{cases}
\displaystyle \sum_{i=0}^{n/2} e^{(\ell)}_i x^{2i} & \mbox{ if $n$ is even, }\\
 \displaystyle \sum_{i=0}^{(n-1)/2} e^{(\ell)}_i x^{2i+1} & \mbox{ if $n$ is odd. }
\end{cases}
\end{align*}
Here $E'_i$ is the sum of all permutations with $i-1$ interior peaks and $E^{(\ell)}_i$ is the sum of all permutations with $i-1$ left peaks. We denote the linear span of the $E'_i$ by $\mathfrak{p}_n$, and the span of the $E^{(\ell)}_i$ by $\mathfrak{p}^{(\ell)}_n$. Theorems \ref{thm:interior} and \ref{thm:left} will establish that the elements $e'_i$, $e_i^{(\ell)}$ defined above are mutually orthogonal idempotents, and with a little more work we will see that they form bases for $\mathfrak{p}_n$ and $\mathfrak{p}^{(\ell)}_n$, respectively.

Let $\eta$ be the involution defined by $\eta(i) = n+1-i$. Then we will also define the \emph{exterior} and \emph{right enriched order polynomials}, denoted $\eom'(P;x)$ and $\rom(P;x)$ respectively, with the relations (shown in section \ref{sec:typeAepp}) that
\begin{align}\label{eq:lomrom}
\Omega'(\pi;x) = \eom'(\eta\pi;x) \mbox{ and } \lom(\pi;x) = \rom(\pi\eta;x).
\end{align}
We can use these polynomials to construct exterior and right structure polynomials $\overline{\rho}(x)$ and $\rho^{(r)}(x)$ in the natural way. We let $\overline{E}'_i$ denote the sum of all permutations with $i$ exterior peaks; $E^{(r)}_i$ denotes the sum of all permutations with $i-1$ right peaks. The span of the $\overline{E}'_i$ is $\overline{\mathfrak{p}}_n$. Though the span of the $E^{(r)}_i$ is not an algebra, we will see that the multiplicative closure of their span is a commutative algebra that contains $\mathfrak{p}^{(\ell)}_n$ as a proper subalgebra. We will show that the coefficients of $\overline{\rho}(x)$, $\overline{e}'_i$, are orthogonal idempotents, and how the $e^{(r)}_i$, coefficients of $\rho^{(r)}(x)$, multiply together and with the $e^{(\ell)}_i$. Notice that \eqref{eq:lomrom} implies that $\rho(x) = \eta \overline{\rho}(x)$ and $\rho^{(\ell)}(x) = \rho^{(r)}(x) \eta$.

The proofs omitted here can be found in section \ref{sec:peak}.

\begin{thm}\label{thm:interior}
As polynomials in $x$ and $y$ with coefficients in the group algebra of the symmetric group we have
\begin{align}
\rho(x)\rho(y) & = \rho(xy), \label{eq:int1}\\
\overline{\rho}(x)\overline{\rho}(y) & = \overline{\rho}(xy),\label{eq:int2}\\
\overline{\rho}(x)\rho(y) & = \overline{\rho}(xy), \label{eq:int3}\\
\rho(x)\overline{\rho}(y) & = \rho(xy). \label{eq:int4}
\end{align}
\end{thm}

Note that equations \eqref{eq:int3} and \eqref{eq:int4} follow from \eqref{eq:int1} and \eqref{eq:int2} upon left multiplication by $\eta$.

\begin{cor}\label{cor:interior}
The algebras $\mathfrak{p}_n$ and $\overline{\mathfrak{p}}_n$ are commutative of dimension $\lfloor \frac{n+1}{2} \rfloor$.
\end{cor}

\begin{proof}
We only give the argument for $\mathfrak{p}_n$. The proof of this and similar corollaries follows the same line of reasoning used in the proof of Corollary \ref{cor:desalg}. Again, the only thing we need to show is that $\spn\{ E'_i\} \subset \spn\{e'_i\}$.

The key fact (see section \ref{sec:typeAepp}) is that for nonnegative integers $k$, $\Omega(\pi;k) = 0$ if $k \leq \pe(\pi)$. So we see that \[ \rho(2) = \Omega'(1;1) E'_1 =
\begin{cases}
\displaystyle \sum_{i=1}^{n/2} e'_i 2^{2i} & \mbox{ if $n$ is even, }\\
\displaystyle \sum_{i=1}^{(n+1)/2} e'_i 2^{2i-1} & \mbox{ if $n$ is odd. }
\end{cases}
\]
We next compute $\rho(4), \rho(6), \ldots$, and the result follows by induction just as in Corollary \ref{cor:desalg}.
\end{proof}

While equation \eqref{eq:int1} (resp. \eqref{eq:int2}) establishes that the $e'_i$, (resp. $\overline{e}'_i$) are mutually orthogonal idempotents, equations \eqref{eq:int3} and \eqref{eq:int4} give $e'_i \overline{e}'_j = \delta_{ij} e'_i$ and $\overline{e}'_i e'_j = \delta_{ij} \overline{e}'_i$. Let $\widehat{\mathfrak{p}}_n$ denote the algebra spanned by both the $E'_i$ and the $\overline{E}'_i$. Then $\mathfrak{p}_n$ and $\overline{\mathfrak{p}}_n$ are left ideals in the (non-commutative) algebra $\widehat{\mathfrak{p}}_n$.

\begin{thm}\label{thm:left}
As polynomials in $x$ and $y$ with coefficients in the group algebra of the symmetric group we have
\begin{align}
\rho^{(\ell)}(x)\rho^{(\ell)}(y) & =  \rho^{(\ell)}(xy),\label{eq:left1}\\
\rho^{(r)}(x)\rho^{(r)}(y) & = \rho^{(\ell)}(xy), \label{eq:left2}\\
\rho^{(\ell)}(x)\rho^{(r)}(y) & = \rho^{(r)}(xy), \label{eq:left3}\\
\rho^{(r)}(x)\rho^{(\ell)}(y) & = \rho^{(r)}(xy).\label{eq:left4}
\end{align}
\end{thm}

Note that equations \eqref{eq:left3} and \eqref{eq:left4} follow from \eqref{eq:left1} and \eqref{eq:left2} upon right multiplication by $\eta$.

\begin{cor}\label{cor:left}
The algebra $\mathfrak{p}^{(\ell)}_n$ is commutative of dimension $\lfloor \frac{n}{2} \rfloor + 1$.
\end{cor}

Equations \eqref{eq:left1} establishes that the $e^{(\ell)}_i$ are mutually orthogonal idempotents. Taken together with \eqref{eq:left2}, \eqref{eq:left3}, and \eqref{eq:left4} we have $e^{(\ell)}_i e^{(\ell)}_j = e^{(r)}_i e^{(r)}_j = \delta_{ij} e^{(\ell)}_i$ and $e^{(\ell)}_i e^{(r)}_j = e^{(r)}_j e^{(\ell)}_i = \delta_{ij} e^{(r)}_i$. If we let $\mathfrak{p}^{(\ell r)}_n$ denote the multiplicative closure of the span of the $E_{i}^{(r)}$, we can see that $\mathfrak{p}^{(\ell)}_n \subsetneq \mathfrak{p}^{(\ell r)}_n$.


\begin{thm}\label{thm:peakideal}
As polynomials in $x$ and $y$ with coefficients in the group algebra of the symmetric group we have
\begin{align}
\rho(x)\rho^{(\ell)}(y) = \rho^{(\ell)}(y)\rho(x) & =  \rho(xy),\label{eq:peakideal1}\\
\overline{\rho}(x)\rho^{(\ell)}(y) = \rho^{(\ell)}(y)\overline{\rho}(x) & =  \overline{\rho}(xy) \label{eq:peakideal2}\\
\rho(x)\rho^{(r)}(y) = \rho^{(r)}(y)\overline{\rho}(x) & = \rho(xy) \label{eq:peakideal3} \\
\overline{\rho}(x)\rho^{(r)}(y) = \rho^{(r)}(y)\rho(x) & = \overline{\rho}(xy) \label{eq:peakideal4}
\end{align}
\end{thm}

Reminiscent of Theorem \ref{thm:idealB}, Theorem \ref{thm:peakideal} (equation \eqref{eq:peakideal1}) implies that $e'_i e^{(\ell)}_j = e^{(\ell)}_j e'_i = \delta_{ij} e'_i$. Let $\dot{\mathfrak{p}}_n$ be the algebra spanned by both $E'_i$ and $E^{(\ell)}_i$. Theorem \ref{thm:peakideal} shows that the peak algebras $\mathfrak{p}_n$, $\mathfrak{p}_n^{(\ell)}$ have the same relationship as the type B descent algebras $\mathfrak{e}^{(c)}_{B,n}$ and $\mathfrak{e}_{B,n}$. Specifically, $\mathfrak{p}_n$, the interior peak algebra, is an ideal in $\dot{\mathfrak{p}}_n$. See \cite{AguiarBergeronNyman} for more connections between type B descent algebras and type A peak algebras. We remark that equation \eqref{eq:peakideal2} implies that the same could be done by replacing $E'_i$ with $\overline{E}'_i$ and $e'_i$ with $\overline{e}'_i$. More conclusions can be drawn, though we will not provide an exhaustive listing of them here. See Table \ref{table:2} for multiplication tables for the various coefficients.

\begin{cor}\label{cor:peakideal}
The algebra $\dot{\mathfrak{p}}_n$ is commutative of dimension $n$ if $n$ is even, dimension $n+1$ if $n$ odd.
\end{cor}

\begin{proof}
We follow the idea from the proof of Corollary \ref{cor:idealB}.
For $i = 1,2,\ldots,\lfloor \frac{n+1}{2} \rfloor$, let $F_{i}^{1}$ be the sum of all permutations with $i-1$ interior peaks and $1 \in \Des(\pi)$, let $F_{i}^{0}$ be the sum of all permutations $\pi$ with $i-1$ interior peaks and $1 \notin \Des(\pi)$. Note that if $n$ is odd, then $F^1_{\frac{n+1}{2}} = 0$. Then
\begin{enumerate}
\item $E^{(\ell)}_1 =  F_{1}^{0}$,
\item $E^{(\ell)}_{\lfloor\frac{n}{2}\rfloor + 1} =
\begin{cases}
F^1_{n/2} & \mbox{ if $n$ is even, }\\
F^1_{\frac{n-1}{2}} + F^0_{\frac{n+1}{2}} & \mbox{ if $n$ is odd, }
\end{cases}
$
\item $E^{(\ell)}_{i} = F_{i-1}^{1} + F_{i}^{0}$ for  $1 < i < \lfloor \frac{n}{2} \rfloor +1$, and
\item $E'_{i} = F_{i}^{1} + F_{i}^{0}$ for $1 \leq i \leq \lfloor \frac{n+1}{2} \rfloor$.
\end{enumerate}
So we see that the $F_{i}^{1}$, $F_{i}^{0}$, which are obviously linearly independent, span the $E^{(\ell)}_{i}$, $E'_{i}$.
\end{proof}

\begin{thm}\label{thm:interiordescent}
As polynomials in $x$ and $y$ with coefficients in the group algebra of the symmetric group we have
\begin{align}
\rho(x)\phi(y) & = \rho(xy), \label{eq:interiordescent}\\
\overline{\rho}(x)\phi(y) &= \overline{\rho}(xy). \label{eq:exteriordescent}
\end{align}
\end{thm}

Theorem \ref{thm:interiordescent} (equation \eqref{eq:interiordescent}) tells us that $e'_i e_j = \delta_{ij} e'_i$. However, the theorem does not imply that the idempotents commute, and in fact the algebra $\mathfrak{ep}_n$ spanned by $E_i$ (the sum of permutations with $i-1$ descents) and $E'_i$ (the sum of permutations with $i-1$ interior peaks) is not commutative for $n >2$. However, the algebra $\mathfrak{p}_n$ is a commutative left ideal in $\mathfrak{ep}_n$. By \eqref{eq:exteriordescent}, the analogous relationship holds for $\overline{\mathfrak{p}}_n$. We note that for $n=3$, the product $\phi(y)\rho(x)$ is not a polynomial in $xy$, and it is not clear if the coefficients of this product have any combinatorial or algebraic significance.

\begin{table}[h]
\begin{tabular}{c|c c c c c}
($\delta_{ij}e^{*}_i e^{*}_j$) & $e_j$ & $e'_j$ & $\overline{e}'_j$ & $e^{(\ell)}_j$ & $e^{(r)}_j$ \\
\hline
$e_i$ & $e_i$ & ? & ? & ? & ? \\

$e'_i$ & $e'_i$ & $e'_i$ & $e'_i$ & $e'_i$ & $e'_i$ \\

$\overline{e}'_i$ & $\overline{e}'_i$ & $\overline{e}'_i$ & $\overline{e}'_i$ & $\overline{e}'_i$ & $\overline{e}'_i$ \\

$e^{(\ell)}_i$ & ? & $e'_i$ & $\overline{e}'_i$ & $e^{(\ell)}_i$ & $e^{(r)}_i$ \\

$e^{(r)}_i$ & ? & $\overline{e}'_i$ & $e'_i$ & $e^{(r)}_i$ & $e^{(\ell)}_i$

\end{tabular}
\vspace{.5cm}
\caption{Multiplication table for type A coefficients. \label{table:2}}
\end{table}

\subsection{Type B}

For the hyperoctahedral group, we will consider only one type of peak. We say a signed permutation $\pi \in \mathfrak{B}_n$ has a peak in position $i = 1,2,\ldots, n-1$ if $\pi(i-1)<\pi(i) > \pi(i+1)$, where, as in our earlier dealings with signed permutations, we require that $\pi(0) = 0$. We will denote the set of peaks by $\Pe_B(\pi)$, and the number of peaks by $\pe_B(\pi)$. For example, the permutation $\pi = (-2,4,-5,3,1)$ has $\Pe_B(\pi) = \{2,4\}$ and $\pe_B(\pi) = 2$. Note that the number of peaks of a signed permutation is between zero and $\lfloor n/2 \rfloor$.

A natural guess at the definition of an Eulerian peak algebra of type B might be the span of sums of permutations with the same number of peaks. However, this definition simply does not work. The following definition does work. Define the elements $E_{i}^{+}, E_{i}^{-}$ in the group algebra of the hyperoctahedral group by:
\begin{align*}
E_{i}^{+} & := \sum_{\substack{\pe(\pi) = i\\ \pi(1) > 0}} \pi, \\
E_{i}^{-} & := \sum_{\substack{\pe(\pi) = i\\ \pi(1) < 0}} \pi.
\end{align*}
In section \ref{sec:peakB} we will show that the linear span of these elements, denoted $\mathfrak{p}_{B,n}$, forms a commutative subalgebra of the group algebra of the hyperoctahedral group. Note that these elements split the collection of permutations with the same number of peaks into two subsets: those that begin with a positive number and those that begin with a negative number. This splitting of cases is similar to splitting left peaks apart from interior peaks. It is not hard to check that $E_{i}^{+}$ and $E_{i}^{-}$ are nonzero for all $0\leq i < \lfloor n/2 \rfloor$. If $n$ is odd, $E_{\frac{n-1}{2}}^{+}$ and $E_{\frac{n-1}{2}}^{-}$ are both nonzero, but if $n$ is even, $E_{n/2}^{+}$ is nonzero while $E_{n/2}^{-} = 0$. In other words, the set $\{ E_{i}^{\pm} \}$ has cardinality $n+1$ for any $n$.

We now move on to describe orthogonal idempotents for the Eulerian peak algebra of the hyperoctahedral group, establishing its commutativity and dimension. As we shall make precise in section \ref{sec:typeBepp}, for any type B poset $P$ there exists a polynomial $\Omega'_B(P;x)$, the \emph{type B enriched order polynomial}, with the following properties. We will see that $\Omega'_B(\pi;(x-1)/4)$ is a degree $n$ polynomial that depends only on the number of peaks of $\pi$ and the sign of $\pi(1)$. By analogy with type A, let
\begin{align*}
\rho_B(x) := \sum_{\pi \in \mathfrak{B}_n} \Omega'_B(\pi;(x-1)/4) \pi & = \sum_{i=0}^{\lfloor n/2 \rfloor} \left( \Omega'_B(i+;(x-1)/4)E_{i}^{+} + \Omega'_B(i-;(x-1)/4)E_{i}^{-} \right) \\
 & = \sum_{i=0}^{n} e'_i x^i,
\end{align*} where $\Omega'_B(i+;x)$ is the enriched order polynomial for any permutation $\pi$ with $i$ peaks and $\pi(1)>0$, $\Omega'_B(i-;x)$ is defined similarly for $\pi$ such that $\pi(1) < 0$. We have the following theorem.

\begin{thm}\label{thm:peakalg2}
As polynomials in $x$ and $y$ with coefficients in the group algebra of the hyperoctahedral group, we have
\[ \rho_B(x)\rho_B(y) = \rho_B(xy). \]
\end{thm}

Then we get a set of $n+1$ orthogonal idempotents since Theorem \ref{thm:peakalg2} gives $e'_i e'_j = \delta_{ij} e'_i$.

\begin{cor}
The algebra $\mathfrak{p}_{B,n}$ is commutative of dimension $n+1$.
\end{cor}

\begin{proof}
This proof is in the spirit of those for Corollaries \ref{cor:desalg} and \ref{cor:interior}. In section \ref{sec:typeBepp} we will be able to make the following observations about type B enriched order polynomials:
\begin{enumerate}
\item $\Omega(0+,0) = 1$,
\item $\Omega(i+,k) = 0$ for $0 < k < i$,
\item $\Omega(i-,k) = 0$ for $0 \leq k \leq i$, and
\item for any integers $i,k$, $\Omega(i+,-(k+1)) = (-1)^n\Omega(i+,k)$ and $\Omega(i-,-k) = (-1)^n \Omega(i-,k)$.
\end{enumerate}
Now we put these observations to use inductively. First, $\rho_B(1) = E_0^+ = e'_0 + e'_1 + \cdots + e'_n$. Next, we consider $\rho_B(-3) = \Omega(0+,-1)E^+_0 + \Omega(0-,-1)E_0^- = e'_0 -3e'_1 + \cdots + (-3)^n e'_n$. The induction proceeds by next computing $\rho_B(5), \rho_B(-7), \rho_B(9)$, etc. Thus we can establish that $\spn\{E_i^{\pm}\} \subset \spn\{ e'_i\}$.
\end{proof}

\section{Enriched $P$-partitions}\label{sec:epp}

In this section we will give the definitions and basic tools needed to use enriched $P$-partitions. In section \ref{sec:typeAepp}, we give two closely related kinds of enriched $P$-partitions for the symmetric group. In section \ref{sec:typeBepp}, we give a definition for enriched $P$-partitions of type B.

\subsection{Type A}\label{sec:typeAepp}

We now introduce the basic theory of enriched $P$-partitions, much of which is due to Stembridge \cite{Stembridge}, and \emph{left} enriched $P$-partitions which, though new, are in the same spirit.

To begin, Stembridge defines $\mathbb{P}'$ to be the set of nonzero integers with the following total order: \[ -1 < 1 < -2 < 2 < -3 < 3 < \cdots\] We will have use for this set, but we view it as a subset of a similar set. Define $\mathbb{P}^{(\ell)}$ to be the integers with the following total order: \[0 < -1 < 1 < -2 < 2 < -3 < 3 < \cdots\] Then $\mathbb{P}'$ is simply the set of all $i \in \mathbb{P}^{(\ell)}$, $i > 0$. In general, for any totally ordered set $X = \{x_1, x_2, \ldots\}$ we define $X^{(\ell)}$ to be the set \[\{x_0, -x_1, x_1, -x_2, x_2, \ldots\},\] with total order \[x_0 < -x_1 < x_1 < -x_2 < x_2 < \cdots \] (so we can think of $X^{(\ell)}$ as two interwoven copies of $X$ along with a zero element) and define $X'$ to be the set $\{ x \in X^{(\ell)} \,|\, x > x_0 \}$. In particular, for any positive integer $k$, $[k]^{(\ell)}$ is the set \[0 < -1 < 1 < -2 < 2 < \cdots < -k < k,\] and $[k]'$ is \[-1 < 1 < -2 < 2 < \cdots < -k < k.\] For any $x_i \in \{x_0\} \cup X$, we say $x_i \geq 0$, or $x_i$ is \emph{nonnegative}. On the other hand, if $i\neq 0$ we say $-x_i < 0$ and $-x_i$ is \emph{negative}. The absolute value removes any minus signs: $|\pm x| = x$ for any $x \in \{x_0\} \cup X$.

For $i$ and $j$ in $X^{(\ell)}$, we write $i\leq^{+} j$ to mean either $i < j$ in $X^{(\ell)}$, or $i = j \geq 0$. Similarly we define $i \leq^{-} j$ to mean either $i < j$ in $X^{(\ell)}$, or $i = j < 0$. For example, on $\mathbb{P}^{(\ell)}$, we have $\{ i \, | \, i \leq^{+} 3 \} = \{\,0, \pm 1, \pm 2, \pm 3 \,\}$, $\{\, i \, | \, i \leq^{-} 3 \,\} = \{\, 0, \pm 1, \pm 2, -3\,\} = \{\, i \,|\, i \leq^{-} -3 \,\}$, $\{ \,i \,| \,0 \leq^{+} i \leq^{+} 2 \,\} = \{\, 0, \pm 1, \pm 2 \,\}$ and $\{\, i \,| \,0 \leq^{-} i \leq^{+} 2 \,\} = \{ \,\pm 1, \pm 2\, \}$.

\begin{defn}[Enriched $P$-partition]\label{def:epp}
An \emph{enriched $P$-partition} (resp. \emph{left enriched $P$-partition}) is an order-preserving map $f: P \to X'$ (resp. $X^{(\ell)}$) such that for all $i <_{P} j$ in $P$,
\begin{enumerate}
\item $f(i) \leq^{+} f(j)$ only if $ i < j$ in $\mathbb{Z}$,
\item $f(i) \leq^{-} f(j)$ only if $ i > j$ in $\mathbb{Z}$.
\end{enumerate}
\end{defn}

It is helpful to remember that enriched $P$-partitions are the nonzero left enriched $P$-partitions. We let $\mathcal{E}(P)$ denote the set of all enriched $P$-partitions; $\mathcal{E}^{(\ell)}(P)$ denotes the set of left enriched $P$-partitions. When $X$ has a finite number of elements, $k$, then the number of (left) enriched $P$-partitions is finite. In this case, define the \emph{enriched order polynomial}, denoted $\Omega'(P;k)$, to be the number of enriched $P$-partitions $f: P \to X'$. The \emph{left enriched order polynomial}, $\lom(P;k)$, is the number of left enriched $P$-partitions $f: P \to X^{(\ell)}$. These enriched order polynomials play the same role in the study of Eulerian peak algebras that ordinary order polynomials play in the study of Eulerian descent algebras.

Just as with ordinary $P$-partitions, we have the fundamental lemma of enriched
P-partitions.

\begin{lem}
\label{lem:FLEPP}
For any poset $P$, the set of all (left) enriched $P$-partitions is the
disjoint union of all (left) enriched $\pi$-partitions for linear
extensions $\pi$ of $P$. Equivalently,
\begin{align*}
\mathcal{E}(P) & = \coprod_{\pi \in \mathcal{L}(P)} \mathcal{E}(\pi),\\
\mathcal{E}^{(\ell)}(P) & = \coprod_{\pi \in \mathcal{L}(P)} \mathcal{E}^{(\ell)}(\pi).
\end{align*}
\end{lem}

The proof of the lemma is identical to the proof of the analogous statement for ordinary $P$-partitions, and the following corollary is immediate.

\begin{cor}
The (left) enriched order polynomial for a poset $P$ is the sum of the (left) enriched order polynomials for all linear extensions of $P$.
\begin{align*}
\Omega'(P;k) & = \sum_{\pi \in \mathcal{L}(P)} \Omega'(\pi;k),\\
\lom(P;k) & = \sum_{\pi \in \mathcal{L}(P)} \lom(\pi;k).
\end{align*}
\end{cor}

Therefore when studying enriched $P$-partitions it is enough to consider the case where $P$ is a permutation. It is easy to describe the set of all enriched $\pi$-partitions in terms of descent sets. For any $\pi \in \mathfrak{S}_n$ we have
\begin{equation}
\label{eq:epp}
\begin{aligned}
\mathcal{E}(\pi) = \{\, f: [n] \to X' & \mid f(\pi(1)) \leq f(\pi(2)) \leq \cdots \leq f(\pi(n)), \\
&  i \notin \Des(\pi) \Rightarrow f(\pi(i)) \leq^{+} f(\pi(i+1))  \\
&  i \in \Des(\pi) \Rightarrow f(\pi(i)) \leq^{-} f(\pi(i+1))  \, \},
\end{aligned}
\end{equation}
and the analogous description for $\mathcal{E}^{(\ell)}(\pi)$ where we replace $X'$ with $X^{(\ell)}$. Notice that \eqref{eq:epp} looks just like the description of $\mathcal{A}(\pi)$, except that we've traded $X$ for $X'$ or $X^{(\ell)}$, $\leq$ is replaced with $\leq^+$, and $\leq^-$ takes the place of $<$.

When we take $X = [k]$, counting the number of solutions to a set of inequalities like \eqref{eq:epp} is not as straightforward as, say, counting the number of integer solutions to a system of ordinary inequalities. The enriched order polynomial is not a simple binomial coefficient as in the case of ordinary $P$-partitions, but it is still possesses some nice properties. We will now present several important properties of enriched order polynomials given by Stembridge \cite{Stembridge}, along with the analogous statements for left enriched order polynomials.

Let $c_l(P)$ denote the number of enriched $P$-partitions $f$ such that $\{\, |f(i)| \,:\, i = 1,2,\ldots,n\,\} = [l]$ as sets. Let $c^0_l(P)$ denote the number of left enriched $P$-partitions $f$ such that $\{\,|f(i)| \,:\, i=1,2,\ldots,n\,\} = [0,l]$. Then we have the following formulas for the enriched order polynomials:
\begin{align*}
\Omega'(P;k) & = \sum_{l=1}^{n} \binom{k}{l} c_{l}(P), \\
\lom(P;k) & = \Omega'(P;k) + \sum_{l=0}^{n-1} \binom{k}{l} c^0_l(P).
\end{align*}
This formula quickly shows that $\Omega'(P;x)$, $\lom(P;x)$ are indeed polynomials, and that they have degree $n$. Though it may not be obvious in this formulation, Stembridge observes (\cite{Stembridge}, Proposition 4.2) that enriched order polynomials satisfy a reciprocity relation. The left enriched order polynomials also satisfy a reciprocity relation, though it requires a shift.
\begin{prp}
\label{prp:eer} We have
\begin{align*}
\Omega'(P;-x) & = (-1)^{n}\Omega'(P;x),\\
\lom(P;-x-1/2) & = (-1)^{n} \lom(P;x-1/2).
\end{align*}
\end{prp}
The proof of this proposition is omitted, though we will say it is straightforward given the generating functions in Theorems \ref{thm:epgf} and \ref{thm:eepgf} below. Before we get too far ahead of the story, we have yet to say why enriched order polynomials are useful for studying peaks of permutations. A hint lies in the fact that by Proposition \ref{prp:eer} we know the number of nonzero terms of $\Omega'(P;x)$ is at most $\lfloor (n+1)/2 \rfloor$, or the number of possible interior peak numbers for permutations in $\mathfrak{S}_n$ (with equality when $P$ is the identity permutation). Likewise, there are $\lfloor n/2 \rfloor +1$ left peak numbers and the left order polynomial $\lom(P;x-1/2)$ has at most $\lfloor n/2 \rfloor + 1$ nonzero terms (it can have a nonzero constant term when $n$ is even).

From \eqref{eq:epp} it is clear that enriched $\pi$-partitions depend on the descent set of $\pi$. A less obvious fact is that they depend only on the set of interior peaks (\cite{Stembridge}, Proposition 2.2). As seen in Theorem \ref{thm:epgf} below, the enriched order polynomial (the \emph{number} of enriched $P$-partitions) depends only on the \emph{number} of interior peaks. Here we give only the generating function for enriched order polynomials of permutations, and remark that by the fundamental Lemma \ref{lem:FLEPP}, we can obtain the order polynomial generating function for any poset by summing the generating functions for its linear extensions.

\begin{thm}[Stembridge \cite{Stembridge}, Theorem 4.1]\label{thm:epgf}
We have the following generating function for enriched order polynomials:
\[ \sum_{k \geq 0} \Omega'(\pi;k)t^k = \frac{1}{2}\frac{(1+t)^{n+1}}{(1-t)^{n+1}}  \left(\frac{4t}{(1+t)^2}\right)^{\pe(\pi)+1}\]
\end{thm}

Notice that this formula implies that $\Omega'(\pi;x)$ has no constant term. We present the proof below since it helps to understand subsequent proofs for left and type B enriched order polynomials.

\begin{proof}
Fix any permutation $\pi \in \mathfrak{S}_n$. From the general theory of $P$-partitions (see, e.g., \cite{Stanley}, chapter 4), we have the following formula for the generating function for ordinary order polynomials: \[\sum_{k \geq 0} \Omega(\pi;k) t^k = \frac{ t^{\des(\pi) + 1}}{(1-t)^{n+1}}\] For any set of integers $D$, let $D+1$ denote the set $\{ d+1 \,\mid\, d \in D\}$. From \cite{Stembridge}, Proposition 3.5, we see that an enriched order polynomial can be written as a sum of ordinary order polynomials: \[ \Omega'(\pi;k) = 2^{\pe(\pi)+1} {\kern -10pt} \sum_{\substack{ D \subset [n-1] \\ \Pe(\pi) \subset D \vartriangle (D + 1)}} {\kern -10pt} \Omega(D;k),\] where $\Omega(D;k)$ denotes the ordinary order polynomial of any permutation with descent set $D$, and $\vartriangle$ denotes the symmetric difference of sets: $A \vartriangle B = (A \cup B)\backslash(A \cap B)$. Putting these two facts together, we get:
\begin{align*}
\sum_{k \geq 0} \Omega'(\pi;k) t^k & = \sum_{k \geq 0} 2^{\pe(\pi)+1} {\kern -10pt} \sum_{\substack{ D \subset [n-1]  \\ \Pe(\pi) \subset D \vartriangle (D + 1)}} {\kern -10pt} \Omega(D;k) t^k \\
 & = 2^{\pe(\pi) +1} {\kern -10pt} \sum_{\substack{ D \subset [n-1]  \\ \Pe(\pi) \subset D \vartriangle (D + 1)}} \sum_{k\geq 0} \Omega(D;k) t^k \\
 & = \frac{ 2^{\pe(\pi)+1} }{ (1-t)^{n+1} } \cdot t {\kern -10pt} \sum_{\substack{ D \subset [n-1] \\ \Pe(\pi) \subset D \vartriangle (D + 1)}} {\kern -10pt} t^{|D|}
\end{align*}
It is not hard to write down the generating function for the sets $D$ by size. We have, for any $j \in \Pe(\pi)$, exactly one of $j$ or $j-1$ is in $D$. There are $n-2\pe(\pi)-1$ remaining elements of $[n-1]$, and they can be included in $D$ or not:
\begin{align*}
\sum_{\substack{ D \subset [n-1] \\ \Pe(\pi) \subset D \vartriangle (D + 1)}} {\kern -10pt} t^{|D|} & =  \underbrace{(t+t)(t+t)\cdots(t+t)}_{\pe(\pi)}\underbrace{(1+t)(1+t)\cdots(1+t)}_{n-2\pe(\pi)-1}\\
& =  (2t)^{\pe(\pi)}(1+t)^{n-2\pe(\pi) - 1}
\end{align*}

Putting everything together, we get \[ \sum_{k \geq 0} \Omega'(\pi;k)t^k = \frac{1}{2}\frac{(1+t)^{n+1}}{(1-t)^{n+1}} \left(\frac{4t}{(1+t)^2}\right)^{\pe(\pi)+1}\] as desired.
\end{proof}

We now derive the generating function for the left enriched polynomials to show they depend only on the number of left peaks. As before, we write down the case where the poset is a permutation.

\begin{thm}\label{thm:eepgf}
We have the following generating function for left enriched order polynomials:
\[ \sum_{k \geq 0} \lom(\pi;k)t^k = \frac{(1+t)^{n}}{(1-t)^{n+1}}\left(\frac{4t}{(1+t)^2}\right)^{\lpe(\pi)}\]
\end{thm}

The proof of Theorem \ref{thm:eepgf} relies on results for quasisymmetric functions from section \ref{sec:qsym}.

\begin{proof}
Fix any permutation $\pi \in \mathfrak{S}_n$. The key fact is given by Theorem \ref{thm:fun} found in section \ref{sec:qsym}: \[ \lom(\pi;k) = 2^{\lpe(\pi)} {\kern -10pt} \sum_{\substack{ D \subset [0,n-1] \\ \lPe(\pi) \subset D \vartriangle (D + 1)}} {\kern -10pt} \Omega_B(D;k),\] where $\Omega_B(D;k)$ denotes the ordinary type B order polynomial of any signed permutation with descent set $D$. It may seem strange to express a type A polynomial related to peaks in terms of type B polynomials related to descents, but as may be more clear later on, left peaks are basically a special case of type B peaks, which are quite naturally related to type B descents.

The generating function for type B order polynomials is (see Reiner \cite{Reiner} for example) \[\sum_{k \geq 0} \Omega_B(\pi;k) t^k = \frac{t^{\des_B(\pi)}}{(1-t)^{n+1}}\] As before, we put these two facts together to get:
\begin{align*}
\sum_{k \geq 0} \lom(\pi;k) t^k & = \frac{2^{\lpe(\pi)}}{(1-t)^{n+1}} {\kern -10pt} \sum_{\substack{ D \subset [0,n-1] \\ \lPe(\pi) \subset D \vartriangle (D + 1)}} {\kern -10pt} t^{|D|} \\
 & = \frac{2^{\lpe(\pi)}}{(1-t)^{n+1}}(2t)^{\lpe(\pi)}(1+t)^{n-2\lpe(\pi)}
\end{align*}

By rearranging terms, we get \[ \sum_{k \geq 0} \lom(\pi;k)t^k = \frac{(1+t)^{n}}{(1-t)^{n+1}}\left(\frac{4t}{(1+t)^2}\right)^{\lpe(\pi)}\] as desired.
\end{proof}

Now we will quickly outline the right and exterior enriched $P$-partitions. We omit unimportant parts of the theory, focusing on the enriched order polynomials, as it is these that are most useful for studying the Eulerian peak algebras. Most of their properties are easily deduced from the properties of the left and interior enriched order polynomials. The only difference is the image set of our enriched $P$-partitions.

The \emph{right enriched order polynomial}, $\rom(P;k)$, is the number of enriched $P$-partitions $f: P \to [k]'\cup \{-(k+1)\}$, and the \emph{exterior enriched order polynomial} is the number of enriched $P$-partitions $f: P \to [k-1]^{(\ell)} \cup \{ -k\}$, denoted $\eom'(P;k)$. Using the line of reasoning mentioned in the introduction, it is not hard to verify that the following proposition is true.

\begin{prp}\label{prp:leftintvsrightext}
We have the following equality of enriched order polynomials,
\begin{align*}
\lom(\pi;x) & =  \rom(\pi\eta;x),\\
\Omega'(\pi;x) & =  \eom'(\eta\pi;x),
\end{align*}
where $\eta$ is the involution defined by $\eta(i) = n+1-i$.
\end{prp}

\begin{proof}
We show that $\Omega'(\pi;x) = \eom'(\eta\pi;x)$ here, and remark that proving the other equality is equally straightforward. Notice that if $(a_1, \ldots, a_n)$ is a solution to
\begin{align*}
0 \leq i_1 \leq \cdots \leq i_n \leq -k, \mbox{ where } s \in \Des(\pi) \Rightarrow i_s \leq^+ i_{s+1}, \\
 s \notin \Des(\pi) \Rightarrow i_s \leq^- i_{s+1},
\end{align*}
then $(a_1', \ldots, a'_n)$ is a solution to
\begin{align*}
-1 \leq i_1 \leq \cdots \leq i_n \leq k, \mbox{ where } s \in \Des(\eta\pi) \Rightarrow i_s \leq^+ i_{s+1}, \\
 s \notin \Des(\eta\pi) \Rightarrow i_s \leq^- i_{s+1},
\end{align*}
and we define $a'$ by $a' = -a-1$ if $a \geq 0$, $a' = -a$ if $a < 0$.
\end{proof}

So the polynomials $\rom(\pi;x)$ and $\eom'(\pi;x)$ inherit all the nice properties outlined for the left and interior enriched order polynomials. In particular, because $\lpe(\pi\eta) = \rpe(\pi)$ and we know that $\lom(\pi\eta;x)$ depends only on the number of left peaks of $\pi\eta$, then it must be that $\rom(\pi;x)$ depends only on the number of right peaks of $\pi$. Because $\pe(\eta\pi) + 1 = \overline{\pe}(\pi)$ and $\Omega'(\eta\pi;x)$ depends only on the number of interior peaks of $\eta\pi$, then $\eom'(\pi;x)$ depends only on the number of exterior peaks of $\pi$.

So while we may not have the enriched order polynomials given by a simple binomial coefficient as with ordinary order polynomials, we do know that we have polynomials that depend only on the number of peaks, and that have as many terms as there are realizable peak numbers. Recall that this is very similar to the case of descents, where we knew that our order polynomials depended on the number of descents, and that the number of terms in these polynomials corresponded to the number of realizable descent numbers. It remains to show that the structure polynomials multiply as stated in the theorems of section \ref{sec:results}. We delay these proofs until section \ref{sec:peak}. First, we present type B enriched $P$-partitions.

\subsection{Type B}\label{sec:typeBepp}

When working with signed permutations, we need to change our notion of a poset slightly. See Chow \cite{Chow}; this definition is a simpler version of the notion due to Reiner \cite{Reiner}.
\begin{defn}
A \emph{type B poset}, or \emph{$\mathfrak{B}_{n}$ poset}, is a poset $P$ whose elements
are $0, \pm 1, \pm 2, \ldots, \pm n$ such that if $i <_{P} j$ then
$-j <_{P} -i$.
\end{defn}

\begin{figure} [h]
\centering
\includegraphics[scale=.8]{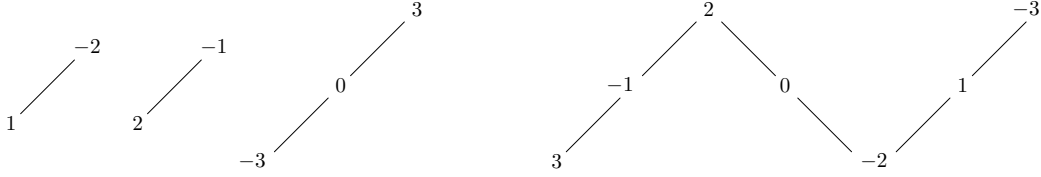}
\caption{Two $\mathfrak{B}_{3}$ posets.\label{fig:typeBposet}}
\end{figure}

\begin{figure} [h]
\centering
\includegraphics[scale=.8]{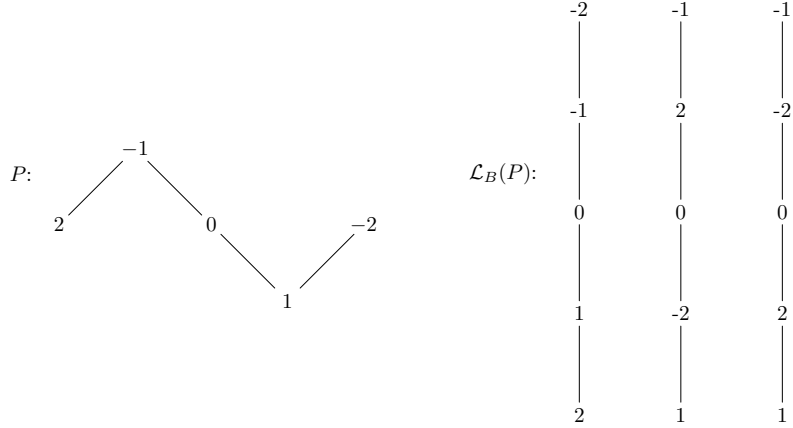}
\caption{A $\mathfrak{B}_{2}$ poset and its linear extensions.\label{fig:typeBlinearextension}}
\end{figure}

Note that if we are given a poset with $n+1$ elements labeled by $0, a_{1},\ldots, a_{n}$ where $a_{i} = i$ or
$-i$, then we can extend it to a $\mathfrak{B}_{n}$ poset of $2n+1$ elements. For example, the type B posets given in Figure \ref{fig:typeBposet} could be specified by the relations $2 <_P -1, 0 <_P 3$ for the poset on the left and $ 0 >_P -2 <_P 1 <_P -3$ for the one on the right. In the same way, any signed permutation $\pi \in \mathfrak{B}_{n}$ is a $\mathfrak{B}_{n}$ poset under the total order $\pi(s) <_{\pi} \pi(s+1)$, $0\leq s\leq n-1$. If $P$ is a type B poset, let $\mathcal{L}_B(P)$ denote the set of linear extensions of $P$ that are themselves type B posets, as in Figure \ref{fig:typeBlinearextension}. Then $\mathcal{L}_B(P)$ is naturally identified with some set of signed permutations.

We will present some alternate notation for the set $X'$ introduced in section \ref{sec:typeAepp}. Let $X = \{x_1, x_2, \ldots\}$ be any totally ordered set. Then we define the totally ordered set $X'$ to be the set \[\{ x_1^{-1}, x_1, x_2^{-1}, x_2, \ldots\},\] with total order \[x_1^{-1} < x_1 < x_2^{-1} < x_2 < \cdots \] We introduce this new notation because we want to avoid confusion in defining the set \[\mathbb{Z}' = \{ \ldots, -2, -2^{-1}, -1, -1^{-1}, 0, 1^{-1}, 1, 2^{-1}, 2, \ldots \},\] with the total order \[ \cdots -2 < -2^{-1} < -1 < -1^{-1} < 0 < 1^{-1} < 1 < 2^{-1} < 2 < \cdots\] In general, if we define $\pm X = \{ \ldots, -x_2, -x_1, x_0, x_1, x_2, \ldots\}$, we have the total order on $\pm X'$ given by \[ \cdots -x_{2} < -x_{2}^{-1} < -x_{1} < -x_{1}^{-1} < x_0 < x_{1}^{-1} < x_{1} < x_{2}^{-1} < x_{2} < \cdots\] We also have the special case for any positive integer $k$, $\pm[k]'$ has total order \[ -k < -k^{-1} < \cdots < -1 < -1^{-1} < 0 < 1^{-1} < 1 < \cdots < k^{-1} < k.\]

For any $x$ in $\pm X'$, let $\varepsilon(x)$ be the exponent on $x$, and let $|x|$ be a map $\pm X' \to X$ that forgets signs and exponents. For example, if $x = -x_{i}^{-1}$, then $\varepsilon(x) = -1 < 0$ and $|x| = x_{i}$, while if $x = x_{i}$, then $\varepsilon(x) = 1 > 0$ and $|x| = x_i$. For $i=0$, we require $\varepsilon(x_0) = 1 > 0$, $|x_0| = x_0$, and $-x_0 = x_0$. Let $x \leq^{+} y$ mean that $x < y$ in $\pm X'$ or $x = y$ and $\varepsilon(x) > 0$. Similarly define $x \leq^{-} y$ to mean that $x < y$ in $\pm X'$ or $x=y$ and $\varepsilon(x) < 0$.

Another way to think of $\mathbb{Z}'$ is as a total ordering of the integer points on the axes in $\mathbb{Z} \times \mathbb{Z}$: \[\cdots (0,-2) < (-1,0) < (0,-1) < (0,0) < (0,1) < (1,0) < (0,2) <\cdots \] In particular, we have $(k,l) < (k',l')$ in $\mathbb{Z}'$ if $k+l < k'+l'$ (in $\mathbb{Z}$), or if $k = l'$ (also in $\mathbb{Z}$). We have $\varepsilon((k,0)) = 1$ for $k\geq 0$, $\varepsilon((0,k)) = -1$ for $k < 0$, and $|(k,l)| = |k+l|$. To negate a point we simply reflect across the perpendicular axis. Note that we could also use this model to understand $\mathbb{P}^{(\ell)}$ from the previous section as all those points $(i,j)$ with $i+j \geq 0$.

\begin{figure} [h]
\centering
\includegraphics[scale = .7]{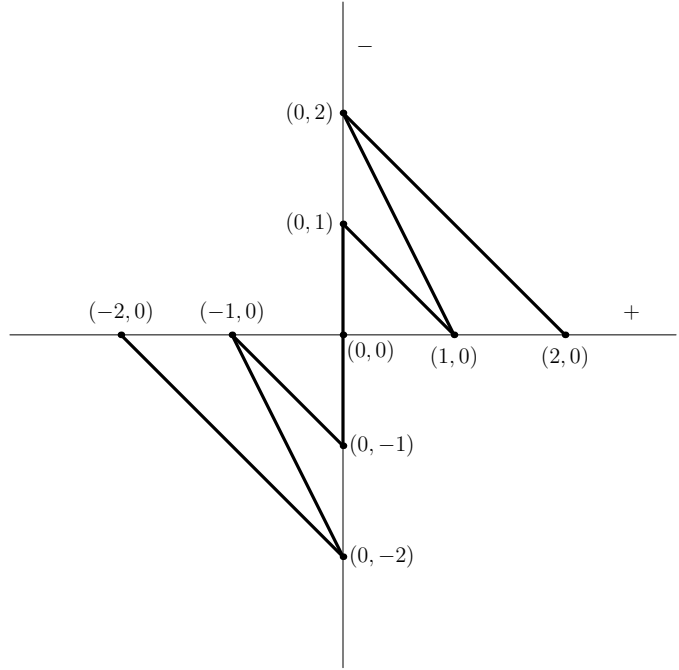}
\caption{One realization of the total order on $\mathbb{Z}'$.}
\end{figure}

\begin{defn}[Type B enriched $P$-partition]
For any $\mathfrak{B}_n$ poset $P$, an \emph{enriched $P$-partition of type B} is an order-preserving map $f : \pm[n] \to \pm X'$ such that for every $i <_{P} j$ in $P$,
\begin{enumerate}
\item $f(i) \leq^{+} f(j)$ only if $i < j$ in $\mathbb{Z}$,
\item $f(i) \leq^{-} f(j)$ only if $i > j$ in $\mathbb{Z}$,
\item $f(-i) = -f(i)$.
\end{enumerate}
\end{defn}

As in the case of ordinary type B $P$-partitions, this definition differs from type A enriched $P$-partitions only in the last condition (see \cite{Petersen}). It forces $f(0) = x_0$, and if we know where to map $a_1, a_2, \ldots, a_n$, where $a_i = i$ or $-i$, then it tells us where to map everything else. In other words, there are $n$ variables. Let $\mathcal{E}_B(P)$ denote the set of all type B enriched $P$-partitions. If we take $X$ to have finite cardinality $k$, then define the \emph{enriched order polynomial of type B}, denoted $\Omega'_B(P;k)$, to be the number of type B enriched $P$-partitions $f: P \to \pm X'$. We have a fundamental lemma.

\begin{lem}
The set of all type B enriched $P$-partitions is the disjoint union of all type B enriched $\pi$-partitions where $\pi$ ranges over all linear extensions of $P$. \[\mathcal{E}_B(P) = \coprod_{\pi \in \mathcal{L}_B(P)} \mathcal{E}_B(\pi).\]
\end{lem}

\begin{cor}
The type B enriched order polynomial for a poset $P$ is the sum of all type B enriched order polynomials for all linear extensions of $P$.
\[\Omega'_B(P;k) = \sum_{\pi \in \mathcal{L}_B(P)} \Omega'_B(\pi;k).\]
\end{cor}

We can easily characterize the type B enriched $\pi$-partitions in terms of descent sets, keeping in mind that if we know where to map $i$, then we know where to map $-i$ by the symmetry property: $f(-i) = -f(i)$. For any signed permutation $\pi \in \mathfrak{B}_n$ we have
\begin{equation}\label{eq:eppsetB}
\begin{aligned}
\mathcal{E}(\pi) = \{\, f: [n] \to \pm X' & \mid  x_0 \leq f(\pi(1)) \leq f(\pi(2)) \leq \cdots \leq f(\pi(n))\\
 &  \quad i \notin \Des_B(\pi) \Rightarrow f(\pi(i)) \leq^{+} f(\pi(i+1)),  \\
 &  \quad i \in \Des_B(\pi) \Rightarrow f(\pi(i)) \leq^{-} f(\pi(i+1)) \,\}.
\end{aligned}
\end{equation}
Notice that since $\varepsilon(x_0) = 1$, then $x_0 \leq^{-} f(\pi(1))$ is the same as saying $x_0 < f(\pi(1))$, and $x_0 \leq^{+} f(\pi(1))$ is the same as $x_0 \leq f(\pi(1))$.

While we have in some sense already said precisely what type B enriched order polynomials are, we need to give a few more properties of them. First of all, let $c_l(P)$ denote the number of type B enriched $P$-partitions $f$ such that $\{\,|f(i)| : i = 1,2,\ldots,n\,\} = [l]$ as sets, and let $c_l^0(P)$ denote the number of type B enriched $P$-partitions $f$ such that $\{\,|f(i)| : i = 1,2,\ldots,n\,\} = [0,l]$. Then we have the following formula for the type B enriched order polynomial: \[\Omega'_B(P;k) = \sum_{l=1}^{n} \binom{k}{l}c_{l}(P) + \sum_{l=0}^{n-1} \binom{k}{l} c_l^0(P).\] This formula shows that the type B enriched order polynomial has degree $n$. If $P = \pi$, a signed permutation with $\pi(1)< 0$, the second term vanishes since $c_{l}^{0}(\pi) = 0$ for all $l$. Notice the similarity between this formula and that of the left order polynomial in the type A case.

We can derive the generating function for type B enriched order polynomials in much the same way as the type A case. From \eqref{eq:eppsetB} it is clear that type B enriched $\pi$-partitions depend only on the descent set of $\pi$. In section \ref{sec:qsym} we will see that they depend on the set of peaks and the sign of $\pi(1)$. Now we will show how the number of type B enriched $\pi$-partitions depends precisely on the number of peaks and the sign of $\pi(1)$. Define $\varsigma(\pi)$ so that $\varsigma(\pi) = 0$ if $\pi(1)$ is positive, $\varsigma(\pi) = 1$ if $\pi(1)$ is negative.

\begin{thm}\label{thm:bepgf}
We have the following generating function for type B enriched order polynomials:
\begin{align}
 \sum_{k \geq 0} \Omega'_B(\pi;k)t^k & =  \frac{(1+t)^{n}}{(1-t)^{n+1}} \left(\frac{2t}{1+t}\right)^{\varsigma(\pi)}  \left(\frac{4t}{(1+t)^2}\right)^{\pe_B(\pi)}\label{eq:bepgf} \\
& =  \left(\frac{1}{2}\right)^{\varsigma(\pi)}  \frac{(1+t)^{n+\varsigma(\pi)}}{(1-t)^{n+1}}  \left(\frac{4t}{(1+t)^2}\right)^{\pe_B(\pi)+\varsigma(\pi)}\nonumber
\end{align}
\end{thm}

Formula \eqref{eq:bepgf} implies that $\Omega'_B(\pi;x)$ depends on \emph{both} the number of peaks and the sign of $\pi(1)$. The similarity between this generating function and the generating functions for type A enriched order polynomials is striking:
\begin{align*}
\mbox{ (Interior peaks) } \quad \sum_{k \geq 0} \Omega'(\pi;k)t^k & =  \frac{1}{2}\frac{(1+t)^{n+1}}{(1-t)^{n+1}}\left(\frac{4t}{(1+t)^2}\right)^{\pe(\pi)+1} \\
\mbox{ (Left peaks) } \quad \sum_{k \geq 0} \lom(\pi;k)t^k  & =  \frac{(1+t)^{n}}{(1-t)^{n+1}}\left(\frac{4t}{(1+t)^2}\right)^{\lpe(\pi)}
\end{align*}

Thus we have the following reciprocity relations, where we recall that $\Omega'_B(i+,x)$ is the enriched order polynomial for any signed permutation with $i$ peaks and $\pi(1) > 0$, $\Omega'_B(i-,x)$ is the enriched order polynomial for any signed permutation with $i$ peaks and $\pi(1) < 0$.

\begin{prp}
We have
\begin{align*}
\Omega'_B(i+;-x-1/2) & =  (-1)^{n}\Omega'_B(i+;x-1/2),\\
\Omega'_B(i-;-x) & =  (-1)^{n} \Omega'_B(i-;x).
\end{align*}
\end{prp}

\begin{proof}[Proof of Theorem \ref{thm:bepgf}]
Fix any permutation $\pi \in \mathfrak{B}_n$. We have the following formula for the generating function of ordinary order polynomials of type B (see Reiner \cite{Reiner}): \[\sum_{k \geq 0} \Omega_B(\pi;k) t^k = \frac{ t^{\des_B(\pi)}}{(1-t)^{n+1}}\] From Theorem \ref{thm:fun} in section \ref{sec:qsym}, we see that \[ \Omega'_B(\pi;k) = 2^{(\pe_B(\pi)+\varsigma(\pi))} {\kern -10pt} \sum_{\substack{ D \subset [0,n-1] \\ \Pe_B(\pi) \subset D \vartriangle (D + 1) \\ \pi(1)< 0 \Rightarrow 0 \in D }} {\kern -10pt} \Omega_B(D;k),\] where $\Omega_B(D;k)$ denotes the ordinary type B order polynomial of any signed permutation with descent set $D$. Putting these two facts together, we get:
\[
\sum_{k \geq 0} \Omega'_B(\pi;k) t^k  =  \frac{ 2^{\pe_B(\pi)+\varsigma(\pi)} }{ (1-t)^{n+1} } \sum_{\substack{ D \subset [0,n-1]  \\ \Pe_B(\pi) \subset D \vartriangle (D + 1) \\ \pi(1)< 0 \Rightarrow 0 \in D }} {\kern -10pt} t^{|D|}
\]
To obtain the generating function for the sets $D$ by size, we proceed in two cases. If we don't require that $0$ is in $D$, that is, if $\pi(1)$ is positive, then we get $(2t)^{\pe_B(\pi)}(1+t)^{n-2\pe_B(\pi)}$ exactly as in the type A case. If $\pi(1) < 0$, we have that $0$ is always in $D$ (and hence $|D| > 0$), while for any $j \in \Pe_B(\pi)$, $j$ must be greater than 1 and exactly one of $j$ or $j-1$ will be in $D$. There are $n-2\pe_B(\pi)-1$ remaining elements of $\{0\}\cup[n-1]$, and they can be included in $D$ or not:
\begin{align*}
\sum_{\substack{ D \subset [0,n-1]  \\ \Pe_B(\pi) \subset D \vartriangle (D + 1) \\ \pi(1)< 0 \Rightarrow 0 \in D }} {\kern -10pt} t^{|D|} & =  t\underbrace{(t+t)(t+t)\cdots(t+t)}_{\pe_B(\pi)}\underbrace{(1+t)(1+t)\cdots(1+t)}_{n-2\pe_B(\pi)-1}\\
& =  t(2t)^{\pe(\pi)}(1+t)^{n-2\pe_B(\pi) - 1}
\end{align*}

Taking the two cases together, we can write \[\sum_{\substack{ D \subset [0,n-1]  \\ \Pe_B(\pi) \subset D \vartriangle (D + 1) \\ \pi(1)< 0 \Rightarrow 0 \in D }} {\kern -10pt} t^{|D|} = t^{\varsigma(\pi)}(2t)^{\pe_B(\pi)}(1+t)^{n-2\pe_B(\pi)-\varsigma(\pi)}\]
Finally, we get \[ \sum_{k \geq 0} \Omega'_B(\pi;k)t^k = \frac{(1+t)^{n}}{(1-t)^{n+1}} \left(\frac{2t}{1+t}\right)^{\varsigma(\pi)} \left(\frac{4t}{(1+t)^2}\right)^{\pe_B(\pi)}\] as desired.
\end{proof}

\subsection{Peak numbers}

Work with Eulerian descent algebras is in a sense a generalization of the study of the Eulerian numbers. Just as there are Eulerian numbers, counting the number of permutations with the same descent number, we also have peak numbers, counting the number of permutations with the same number of peaks (not to be confused with $\pe(\pi)$, the peak number of a permutation). We will not devote much time to this topic, but state only those properties that are easy observations given the theory of enriched $P$-partitions developed in this paper.

We denote the number of permutations in $\mathfrak{S}_n$ with descent number $k$ by the Eulerian number $A_{n,k+1}$, and we recall that the Eulerian polynomial is defined as \[A_{n}(t) = \sum_{\pi\in\mathfrak{S}_{n}} t^{\des(\pi)+1} =
\sum_{i=1}^{n} A_{n,i}t^{i}.\] Similarly, we denote the number of signed permutations in $\mathfrak{B}_n$ with $k$ descents by $B_{n,k}$ and define the \emph{type B Eulerian polynomial} (see \cite{Brenti}) as \[B_n(t) = \sum_{\pi \in \mathfrak{B}_n } t^{\des_B(\pi)}.\] The number of signed permutations with $k$ cyclic descents is $B_{n,k}^{(c)}$ the \emph{type B cyclic Eulerian polynomial} is \[B_{n}^{(c)}(t) = \sum_{\pi\in\mathfrak{B}_{n}} t^{\cdes_B(\pi)} = \sum_{i=1}^{n}B_{n,i}^{(c)}t^{i}.\] We have the following relationship between the Eulerian polynomial and the type B cyclic Eulerian polynomial, proved in \cite{Petersen} using $P$-partitions, and also by Fulman \cite{Fulman}.

\begin{prp}[Petersen \cite{Petersen}]\label{prp:augeul}
The number of signed permutations with $i+1$ cyclic descents is
$2^{n}$ times the number of unsigned permutations with $i$
descents, $0\leq i \leq n-1$. In other words, $B_{n}^{(c)}(t) = 2^{n}A_{n}(t)$.
\end{prp}

We will make some similar observations. Following Stembridge \cite{Stembridge}, we denote the number of permutations of $n$ with $k$ interior peaks by $P_{n,k}$. We define the interior peak polynomial as \[W_{n}(t) = \sum_{\pi\in\mathfrak{S}_{n}} t^{\pe(\pi)+1} = \sum_{i=1}^{\lfloor \frac{n+1}{2} \rfloor} P_{n,i}t^{i}.\] Similarly, we define the left peak polynomial as \[W^{(\ell)}_{n}(t) = \sum_{\pi \in \mathfrak{S}_n} t^{\lpe(\pi)} = \sum_{i=0}^{\lfloor \frac{n}{2} \rfloor} P^{(\ell)}_{n,i}t^{i}.\]

Enriched $P$-partitions give us the tools to prove the following propositions relating peak polynomials to Eulerian polynomials. Proposition \ref{prp:peeul1} appears in Remark 4.8 of \cite{Stembridge}; the second equality follows from Proposition \ref{prp:augeul}.
\begin{prp}[Stembridge \cite{Stembridge}]\label{prp:peeul1}
We have the following relation between the interior peak polynomial, the Eulerian polynomial, and the type B cyclic Eulerian polynomial:
\[ W_{n}\left( \frac{4t}{(1+t)^2} \right) = \frac{2^{n+1}}{(1+t)^{n+1}} A_{n}(t) = \frac{2}{(1+t)^{n+1}} B^{(c)}_{n}(t).\]
\end{prp}
\begin{prp}\label{prp:peeul2}
We have the following relation between the left peak polynomial, and the type B Eulerian polynomial:
\[
 W^{(\ell)}_{n}\left( \frac{4t}{(1+t)^2} \right)  =  \frac{1}{(1+t)^{n}} B_n(t).
\]
\end{prp}

\begin{proof}[Proof of Proposition \ref{prp:peeul1}.]
Recall that we have the following formula for the ordinary Eulerian polynomials (see \cite{Stanley}):
\[\sum_{k\geq 0}k^{n}t^{k} = \frac{A_{n}(t)}{(1-t)^{n+1}}.\]

Now let $P$ be an antichain of $n$ elements labeled $1, 2,\ldots,n$. The number of enriched
$P$-partitions $f: [n] \to [k]'$ is $(2k)^n$ since there are $2k$ elements in $[k]'$ and there are no relations among the elements of the antichain.  Therefore $\Omega'(P;k) =
(2k)^{n}$, and since we have $\mathcal{L}(P) = \mathfrak{S}_{n}$, Theorem \ref{thm:epgf} gives
\[\frac{1}{2} \frac{(1+t)^{n+1}}{(1-t)^{n+1}} W_{n}\left(\frac{4t}{(1+t)^2} \right) = \sum_{k\geq 0}(2k)^{n}t^{k}
= 2^{n}\sum_{k\geq 0}k^{n}t^{k} = \frac{2^{n}A_{n}(t)}{(1-t)^{n+1}}.\] Rearranging terms gives the desired result: \[W_{n}\left( \frac{4t}{(1+t)^2} \right) = \frac{2^{n+1}}{(1+t)^{n+1}} A_{n}(t).\]
\end{proof}

\begin{proof}[Proof of Proposition \ref{prp:peeul2}.]
If we let $P$ be an antichain of $n$ elements, the number of left enriched
$P$-partitions $f: [n] \to [k]^{(\ell)}$ is $(2k+1)^n$ since there are $2k+1$ elements in $[k]^{(\ell)}$ and there are no relations among the elements of the antichain.  Therefore $\lom(P;k) = (2k+1)^{n}$. But if $P_B$ is the antichain on $\pm[n]$, the order polynomial $\Omega_B(P_B;k)$ is also $(2k+1)^n$. Since we have $\mathcal{L}(P) = \mathfrak{S}_{n}$, $\mathcal{L}_B(P_B) = \mathfrak{B}_n$, Theorem \ref{thm:eepgf} gives
\[
\frac{(1+t)^{n}}{(1-t)^{n+1}} W^{(\ell)}_{n}\left(\frac{4t}{(1+t)^2} \right) = \sum_{k\geq 0}(2k+1)^{n}t^{k}
 =  \frac{B_n(t)}{(1-t)^{n+1}}.\\
\]
Rearranging terms gives the desired result:
\[W^{(\ell)}_{n}\left( \frac{4t}{(1+t)^2} \right) = \frac{1}{(1+t)^{n}} B_n(t).\]
\end{proof}

We can define type B peak numbers and type B peak polynomials. We will denote the number of signed permutations of $n$ with $k$ peaks and $\pi(1)>0$ by $P^{+}_{n,k}$. We denote the number of signed permutations of $n$ with $k$ peaks and $\pi(1) < 0$ by $P^{-}_{n,k+1}$. We define the type B peak polynomials by
\begin{align*}
W^{+}_{n}(t) & = \sum_{\substack{\pi\in\mathfrak{B}_{n} \\ \pi(1)>0 }} t^{\pe(\pi)} =  \sum_{i=0}^{\lfloor \frac{n}{2} \rfloor} P^{+}_{n,i}t^{i}\\
W^{-}_{n}(t) & = \sum_{\substack{ \pi \in \mathfrak{B}_n \\ \pi(1)<0 }} t^{\pe(\pi)+1} =   \sum_{i=1}^{\lfloor \frac{n}{2}\rfloor +1} P^{-}_{n,i}t^{i}.
\end{align*}

Similarly to the type A case, we have the following proposition.
\begin{prp}\label{prp:bpeeul1}
We have the following relation between type B peak polynomials and the type B Eulerian polynomial:
\begin{align*}
W^{+}_{n}\left( \frac{4t}{(1+t)^2} \right) &+  \frac{1+t}{2}  W^{-}_{n}\left( \frac{4t}{(1+t)^2} \right)\\ & =  \frac{1}{2(1+t)^n}\left( B_n(\sqrt{t})(1+\sqrt{t})^n + B_n(-\sqrt{t})(1-\sqrt{t})^n \right).
\end{align*}
\end{prp}

\begin{proof}
If we let $P$ be the type B antichain on $\pm[n]$, then we have $\Omega'_B(P;k) = (4k+1)^n = \Omega_B(P;2k)$. Thus,
\begin{align*}
 \frac{(1+t)^{n}}{(1-t)^{n+1}} W^{+}_{n}\left(\frac{4t}{(1+t)^2} \right) & + \frac{1}{2}\frac{(1+t)^{n+1}}{(1-t)^{n+1}} W^{-}_{n}\left(\frac{4t}{(1+t)^2} \right) \\ &  = \sum_{k \geq 0} (4k+1)^n t^k  =  \sum_{k\geq 0} \Omega_B(P;2k) t^k.
\end{align*}
We know that the generating function for $\Omega_B(P;k)$ is given by $F(t) = B_n(t)/(1-t)^{n+1}$. However, we want the generating function for only the even terms. We can get this power series by computing \[ \frac{ F(\sqrt{t}) + F(-\sqrt{t})}{2} \] and the proposition follows by rearranging terms.
\end{proof}

The proofs of these propositions lead naturally to a general proposition, not directly related to peaks. Let $P_i$ denote the type B poset given by the relations $0 <_{P_i} -j$ for $0 \leq j \leq i$, $0 <_{P_i} j$ for $i < j \leq n$. Let
\begin{align*}
 W_{n,i}(t) & =  \binom{n}{i} \sum_{\pi \in \mathcal{L}_B(P_i)} t^{\des_B(\pi)},\\
 & =  \sum_{\pi \in \mathfrak{B}_{n,i}} t^{\des_B(\pi)},
\end{align*}
where $\mathfrak{B}_{n,i}$ denotes the set of all signed permutations with exactly $i$ minus signs. It is then straightforward to prove the following.

\begin{prp}\label{prp:bpeeul2}
Let $\alpha$ be an indeterminate. We have the following:
\[ \frac{1}{(1-t)^{n+1}}\sum_{i=0}^{n} \alpha^i W_{n,i}(t) = \sum_{k\geq 0} ((\alpha+1)k + 1)^n t^k.\]
\end{prp}

\begin{proof}
We can write $((\alpha+1)k + 1)^n = (\alpha k + k + 1)^n$ as \[ \sum_{i=0}^n \binom{n}{i} \alpha^i k^i (k+1)^{n-i}.\]
It is easy to verify that $\Omega_B(P_i; k) = k^i(k+1)^{n-i}$.
\end{proof}

We can interpret these generating functions as type B generating functions weighted by the number of minus signs. So the cases already seen correspond to $\alpha = 1,3$.

\section{Proofs}\label{sec:peak}

We now present the application of the theory of enriched $P$-partitions to the study of commutative peak algebras.

\subsection{Type A structure formulas}\label{sec:interiorandleft}

In this section we will prove Theorems \ref{thm:interior} through \ref{thm:interiordescent}, as described in section \ref{sec:results}. The proofs follow the same basic structure of Theorem \ref{thm:ges}, and we will present fewer and fewer details as we go, focusing on only the crucial differences. For Theorem \ref{thm:interior}, we need to show
\begin{align}
\rho(x)\rho(y) & = \rho(xy),\label{eq:peakform} \\
\overline{\rho}(x)\overline{\rho}(y) & = \overline{\rho}(xy).\label{eq:peakform2}
\end{align}
where $\rho(x)$, $\overline{\rho}(x)$ are as defined in section \ref{sec:results}.

\begin{proof}[Proof of Theorem \ref{thm:interior}]
We will prove \eqref{eq:peakform}. By equating the coefficient of $\pi$ on both sides of equation \eqref{eq:peakform} we know that we need only prove the following claim: For any permutation $\pi \in \mathfrak{S}_n$ and positive integers $k$,$l$, we have \[\Omega'(\pi;2kl) = \sum_{\sigma\tau = \pi} \Omega'(\sigma;k)\Omega'(\tau;l).\]

We will interpret the left-hand side of the equation in such a way that we can split it apart to form the right hand side. Rather than considering $\Omega'(\pi;2kl)$ to count maps $f: \pi \to [2kl]'$, we will understand it to count maps $f: \pi \to [l]'\times[k]'$, where we take the \emph{up-down} order on $[l]'\times[k]'$. The up-down order is defined as follows (see Figure \ref{fig:updown}): $(i,j) < (i',j')$ if and only if
\begin{enumerate}
\item $i < i'$, or

\item $i = i' > 0$ and $j < j'$, or

\item $i = i' < 0$ and $j > j'$.
\end{enumerate}
So if the horizontal coordinate is negative, we read the columns from the top down, if the horizontal coordinate is positive, we read from the bottom up. Then $\Omega'(\pi;2kl)$ is the number of solutions to
\begin{equation}\label{eq:pm}
(-1,k) \leq (i_1,j_1) \leq (i_2,j_2) \leq\cdots \leq (i_n,j_n) \leq (l,k)
\end{equation}
where $(i_s,j_s) \leq^{-} (i_{s+1},j_{s+1})$ if $s \in \Des(\pi)$ and $(i_s,j_s) \leq^{+} (i_{s+1},j_{s+1})$ otherwise. For example, if $\pi = (1,3,2)$, we will count the number of points \[(-1,k) \leq (i_1,j_1) \leq^{+} (i_2,j_2) \leq^{-} (i_3,j_3) \leq (l,k).\] Here we write $(i,j) \leq^{+} (i',j')$ in one of three cases: if $i<i'$, or if $i = i' > 0$ and $j \leq^{+} j'$, or if $i = i' < 0$ and $j \geq^{-} j'$. Similarly, $(i,j) \leq^{-} (i',j')$ if $i < i'$, or if $i=i' > 0$ and $j \leq^{-} j'$, or if $i = i' < 0$ and $j \geq^{+} j'$.

\begin{figure} [h]
\centering
\includegraphics[scale = .7]{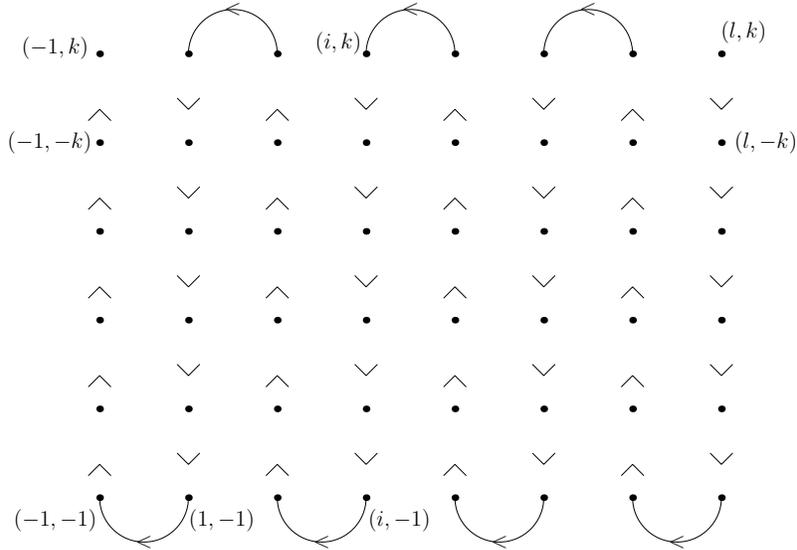}
\caption{The up-down order for $[l]' \times [k]'$.\label{fig:updown}}
\end{figure}

To get the result we desire, we will sort the set of all solutions to \eqref{eq:pm} into distinct cases indexed by subsets $I \subset [n-1]$. The sorting depends on $\pi$ and proceeds as follows. Let $F = ( (i_1, j_1), \ldots, (i_n, j_n) )$ be any solution to \eqref{eq:pm}. For any $s=1,2,\ldots,n-1$, if $\pi(s) < \pi(s+1)$, then $(i_{s},j_{s}) \leq^{+} (i_{s+1},j_{s+1})$, which falls into one of
two mutually exclusive cases:
\begin{align}
i_{s} \leq^{+} i_{s+1} & \mbox{ and }  j_{s}\leq^{+} j_{s+1}, \mbox{ or} \label{eqn:e1}\\
i_{s} \leq^{-} i_{s+1} & \mbox{ and }  j_{s} \geq^{-} j_{s+1}. \label{eqn:e2}
\end{align}
If $\pi(s) > \pi(s+1)$, then $(i_{s},j_{s}) \leq^{-} (i_{s+1},j_{s+1})$,
which we split as:
\begin{align}
i_{s} \leq^{+} i_{s+1} & \mbox{ and }  j_{s} \leq^{-} j_{s+1}, \mbox{ or} \label{eqn:e3}\\
i_{s} \leq^{-} i_{s+1} & \mbox{ and }  j_{s} \geq^{+} j_{s+1},\label{eqn:e4}
\end{align}
also mutually exclusive. Define $I_F$ to be the set of all $s$ such that either \eqref{eqn:e2} or \eqref{eqn:e4} holds for $F$. Notice that in both cases, $i_s \leq^{-} i_{s+1}$. Now for any $I \subset [n-1]$, let $S_I$ be the set of all solutions $F$ to \eqref{eq:pm} satisfying $I_F = I$. We have split the solutions of \eqref{eq:pm} into $2^{n-1}$ distinct cases indexed by all the different subsets $I$ of $[n-1]$.

For any particular $I\subset [n-1]$, form the poset $P_{I}$ of the elements
$1,2,\ldots,n$ by $\pi(s) <_{P_{I}} \pi(s+1)$ if $s \notin I$,
$\pi(s) >_{P_{I}} \pi(s+1)$ if $s\in I$. We form a zig-zag poset (see Figure \ref{fig:zigzag}) of $n$ elements labeled consecutively by $\pi(1),
\pi(2),\ldots,\pi(n)$ with downward zigs corresponding to the
elements of $I$.

For any solution $F$ in $S_I$, let $f: [n] \to [k]'$ be defined by $f(\pi(s)) = j_s$. We will show that $f$ is an enriched $P_{I}$-partition. If $\pi(s) <_{P_{I}} \pi(s+1)$ and $\pi(s) < \pi(s+1)$
in $\mathbb{Z}$, then \eqref{eqn:e1} tells us that $f(\pi(s)) =
j_{s} \leq^{+} j_{s+1} = f(\pi(s+1))$. If $\pi(s) <_{P_{I}} \pi(s+1)$
and $\pi(s) > \pi(s+1)$ in $\mathbb{Z}$, then \eqref{eqn:e3} tells
us that $f(\pi(s)) = j_{s} \leq^{-} j_{s+1} = f(\pi(s+1))$. If $\pi(s)
>_{P_{I}} \pi(s+1)$ and $\pi(s) < \pi(s+1)$ in $\mathbb{Z}$,
then \eqref{eqn:e2} gives us that $f(\pi(s))= j_{s} \geq^{-} j_{s+1} =
f(\pi(s+1))$. If $\pi(s) >_{P_{I}} \pi(s+1)$ and $\pi(s) >
\pi(s+1)$ in $\mathbb{Z}$, then \eqref{eqn:e4} gives us that
$f(\pi(s)) = j_{s} \geq^{+} j_{s+1} = f(\pi(s+1))$. In other words, we
have verified that $f$ is a $P_{I}$-partition. So for any particular solution in $S_I$, the $j_{s}$'s can be thought of as an enriched $P_{I}$-partition.

Conversely, any enriched $P_{I}$-partition $f$ gives a solution in $S_I$ since if $j_s = f(\pi(s))$, then \[(( i_1,j_1),\ldots,(i_n,j_n)) \in S_I\] if and only if $1 \leq i_1 \leq \cdots \leq i_n \leq l$ and $i_s \leq^{-} i_{s+1}$ for all $s \in I$, $i_s \leq^+ i_{s+1}$ for $s \notin I$. We can therefore turn our attention to counting enriched $P_{I}$-partitions.

The remainder of the argument follows the proof of Theorem \ref{thm:ges}.

Equation \eqref{eq:peakform2} is proved in exactly the same fashion. The crucial first step is to understand $\eom'(\pi;2kl)$ as counting enriched $\pi$-partitions $f: \pi \to ([l-1]^{(\ell)}\cup\{-l\}) \times ([k-1]^{(\ell)}\cup \{-k\})$ with the up-down order.
\end{proof}

Now we will prove Theorem \ref{thm:left}. This proof is nearly identical to that of Theorem \ref{thm:interior}---we will highlight only the important differences. We wish to prove the following formulas:
\begin{align}
\rho^{(\ell)}(x)\rho^{(\ell)}(y) & =  \rho^{(\ell)}(xy),\label{eq:leftpeakform} \\
\rho^{(r)}(x)\rho^{(r)}(y) & =  \rho^{(\ell)}(xy).\label{eq:leftpeakform2}
\end{align}
where $\rho^{(\ell)}$, $\rho^{(r)}$ are defined in section \ref{sec:results}.

\begin{proof}[Proof of Theorem \ref{thm:left}]
By equating the coefficient of $\pi$ on both sides of equation \eqref{eq:leftpeakform} we know that we need only prove \[\lom(\pi;2kl+k+l) = \sum_{\sigma\tau = \pi} \lom(\sigma;k)\lom(\tau;l).\]
We will think of the left-hand side of the equation as counting maps $f: \pi \to [l]^{(\ell)}\times[k]^{(\ell)}$, where, as in the proof of Theorem \ref{thm:interior}, we take the up-down order on $[l]^{(\ell)}\times[k]^{(\ell)}$.

Then $\lom(\pi;2kl+k+l)$ is the number of solutions to
\[
(0,0) \leq (i_1,j_1) \leq (i_2,j_2) \leq\cdots \leq (i_n,j_n) \leq (l,k)
\]
where $(i_s,j_s) \leq^{-} (i_{s+1},j_{s+1})$ if $s \in \Des(\pi)$ and $(i_s,j_s) \leq^{+} (i_{s+1},j_{s+1})$ otherwise. Recall that in the up-down order we write $(i,j) \leq^{+} (i',j')$ in one of three cases: if $i<i'$, or if $i = i' \geq 0$ and $j \leq^{+} j'$, or if $i = i' < 0$ and $j \geq^{-} j'$. Similarly, $(i,j) \leq^{-} (i',j')$ if $i < i'$, or if $i=i' \geq 0$ and $j \leq^{-} j'$, or if $i = i' < 0$ and $j \geq^{+} j'$. See Figure \ref{fig:downup}.

\begin{figure} [h]
\centering
\includegraphics[scale = .7]{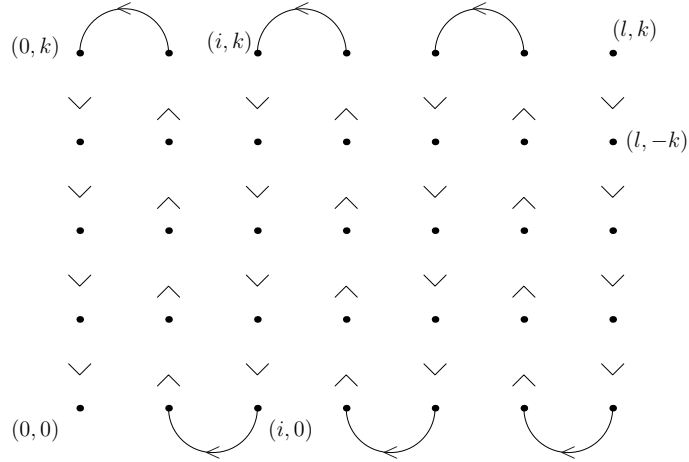}
\caption{The up-down order for $[l]^{(\ell)} \times [k]^{(\ell)}$.\label{fig:downup}}
\end{figure}

The rest of the proof is identical to that of Theorem \ref{thm:interior}.

The proof for equation \eqref{eq:leftpeakform2} follows by considering $\lom(\pi;2kl+k+l)$ as counting enriched $\pi$-partitions $f: \pi \to ([l]'\cup \{-(l+1)\}) \times ([k]'\cup\{-(k+1)\})$ with the up-down order.
\end{proof}

We now give the proof of Theorem \ref{thm:peakideal}. By manipulating the identities $\rho(x) = \eta \overline{\rho}(x)$ and $\rho^{(\ell)}(x) = \rho^{(r)}(x)\eta$, we can see that it will suffice to prove the following identities:
\begin{align}
\rho(x)\rho^{(\ell)}(y) = \rho^{(\ell)}(y)\rho(x) & =  \rho(xy),\label{eq:peakideala}\\
\rho(x)\rho^{(r)}(y) & =  \rho(xy), \label{eq:peakidealb} \\
\rho^{(r)}(y)\rho(x) & = \overline{\rho}(xy). \label{eq:peakidealc}
\end{align}

\begin{proof}[Proof of Theorem \ref{thm:peakideal}]
Conceptually, this proof is little different from the proofs of Theorems \ref{thm:interior} or \ref{thm:left}. We first outline the argument for \eqref{eq:peakideala}. Equating coefficients, we need only prove
\begin{align}
\Omega'(\pi;2kl+l) & = \sum_{\sigma\tau = \pi} \lom(\sigma;k)\Omega'(\tau;l), \label{eq:l'} \\
 & = \sum_{\sigma\tau = \pi} \Omega'(\sigma;l)\lom(\tau;k). \label{eq:'l}
\end{align}
For equation \eqref{eq:l'}, the key is to think of the left-hand side of the equation as counting maps $f: \pi \to [l]' \times[k]^{(\ell)}$, with the up-down order on $[l]' \times[k]^{(\ell)}$.

For \eqref{eq:'l}, we count enriched $\pi$-partitions $f: \pi \to [k]^{(\ell)} \times [l]'$ with the up-down order, and \eqref{eq:peakideala} follows.

To prove \eqref{eq:peakidealb}, we interpret $\Omega'(\pi;2kl+l)$ as counting $\pi$-partitions $f: \pi \to ([l]'\cup \{-(l+1)\}) \times [k]'$ with the up-down order. For \eqref{eq:peakidealc}, we take the up-down order on $[k]' \times ([l]'\cup \{-(l+1)\})$.
\end{proof}

For theorem \ref{thm:interiordescent} we need to show
\begin{align}
\rho(x)\phi(y) = \rho(xy),\label{eq:peakdes}
\end{align}
where $\phi(x)$ is as defined in section \ref{sec:des}.

\begin{proof}[Proof of Theorem \ref{thm:interiordescent}]
This proof mixes the ideas from the proofs of Theorem \ref{thm:ges} and Theorem \ref{thm:interior}. We will equate coefficients of \eqref{eq:peakdes} and show that \[\Omega'(\pi;2kl) = \sum_{\sigma\tau = \pi} \Omega'(\sigma;k) \Omega(\tau;2l)\] for any $\pi \in \mathfrak{S}_n$. First, we consider $\Omega'(\pi;2kl)$ as counting all enriched $P$-partitions $f: [n] \to [l]'\times[k]'$, where now we take the \emph{lexicographic} ordering on the image set (so $(i,j) < (i',j')$ whenever $i < i'$ or $i=i'$ and $j<j'$). An enriched $P$-partition is any solution $((i_1, j_1), \ldots, (i_n, j_n))$ to \[ (-1,-1) \leq (i_1,j_1) \leq (i_2,j_2) \leq \cdots \leq (i_n, j_n) \leq (l,k),\] where $(i_s, j_s) \leq^{+} (i_{s+1},j_{s+1})$ if $s \notin \Des(\pi)$, $(i_s, j_s) \leq^{-} (i_{s+1}, j_{s+1})$ if $s\in \Des(\pi)$. We split all such solutions into cases as before. If $\pi(s) < \pi(s+1)$, then $(i_{s},j_{s}) \leq^{+} (i_{s+1},j_{s+1})$, which falls into one of two mutually exclusive cases:
\begin{align*}
i_{s} \leq i_{s+1} & \mbox{ and } j_{s}\leq^{+} j_{s+1}, \mbox{ or }\\
i_{s} < i_{s+1} & \mbox{ and }  j_{s} \geq^{-} j_{s+1}.
\end{align*}
If $\pi(s) > \pi(s+1)$, then $(i_{s},j_{s}) \leq^{-} (i_{s+1},j_{s+1})$,
which means either:
\begin{align*}
i_{s} \leq i_{s+1} & \mbox{ and }  j_{s} \leq^{-} j_{s+1}, \mbox{ or }\\
i_{s} < i_{s+1} & \mbox{ and }  j_{s} \geq^{+} j_{s+1}.
\end{align*}
With this splitting, the result follows as in the proof of Theorem \ref{thm:ges}. The $i_s$'s are counted with ordinary order polynomials (noting that $[l]'$ has $2l$ elements); enriched order polynomials count the $j_s$'s.
\end{proof}

\subsection{The type B peak structure formula}\label{sec:peakB}

In this section we prove Theorem \ref{thm:peakalg2}. We want to show
\begin{align}
\rho_B(x)\rho_B(y) = \rho_B(xy),\label{eq:peakformb}
\end{align}
where $\rho_B(x)$ is defined in the introduction.

\begin{proof}[Proof of Theorem \ref{thm:peakalg2}.]
This proof is nearly identical to the proofs of the analogous Theorems \ref{thm:interior} and \ref{thm:left}. By equating the coefficient of $\pi$ on both sides of equation \eqref{eq:peakformb} it suffices to prove that for any permutation $\pi \in \mathfrak{B_n}$ and positive integers $k,l$, we have \[\Omega'_B(\pi;4kl+k+l) = \sum_{\sigma\tau = \pi} \Omega'_B(\sigma;k)\Omega'_B(\tau;l).\]

We will interpret $\Omega'_B(\pi;4kl+k+l)$ as counting maps $f: \pi \to \pm[l]' \times \pm[k]'$, where we take the up-down order on $\pm[l]'\times \pm[k]'$. We count up the columns that have positive exponent and down columns with negative exponent. Notice that we can restrict our attention to all the points greater than or equal to $(0,0)$, since everything else is determined by the symmetry property of type B enriched $P$-partitions: $f(-i) = -f(i)$. We consider $\Omega'_B(\pi;4kl+k+l)$ to be the number of solutions to
\begin{equation}\label{eq:pmb}
(0,0) \leq (i_1,j_1) \leq (i_2,j_2) \leq \cdots \leq (i_n,j_n) \leq (l,k)
\end{equation}
where $(i_s,j_s) \leq^{-} (i_{s+1},j_{s+1})$ if $s \in \Des(\pi)$ and $(i_s,j_s) \leq^{+} (i_{s+1},j_{s+1})$ otherwise. For example, if $\pi = (-3,1,-2)$, we will count the number of points $((i_1,j_1), (i_2, j_2), (i_3, j_3))$ such that \[(0,0) \leq^{-} (i_1,j_1) \leq^{+} (i_2,j_2) \leq^{-} (i_3,j_3) \leq (l,k).\] Here $(i,j) \leq^{+} (i',j')$ means  $i<i'$, or if $i = i'$ with $\varepsilon(i) > 0$ and $j \leq^{+} j'$, or if $i = i'$ with $\varepsilon(i) < 0$ and $j \geq^{-} j'$. Similarly, $(i,j) \leq^{-} (i',j')$ if $i < i'$, or if $i=i'$ with $\varepsilon(i) > 0$ and $j \leq^{-} j'$, or if $i = i'$ with $\varepsilon(i) < 0$ and $j \geq^{+} j'$.

\begin{figure} [h]
\centering
\includegraphics[scale = .6]{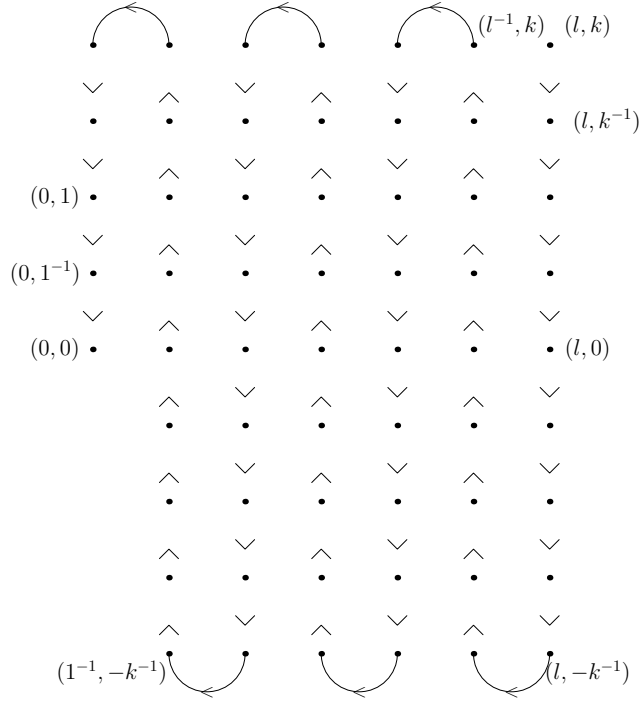}
\caption{The up-down order on $\pm[l]'\times\pm[k]'$ with points greater than or equal to $(0,0)$.}
\end{figure}

Just as with the type A case, we will want to group the solutions to \eqref{eq:pmb} into cases that we will count using enriched order polynomials. Here there are are $2^n$ cases, indexed by subsets of $[0,n-1]$. The grouping depends on $\pi$ and proceeds as follows. Let $F = ( (i_1, j_1), \ldots, (i_n, j_n) )$ be any solution to \eqref{eq:pmb}, and fix $\pi(0) = i_0 = j_0 = 0$. For any $s=0,1,2,\ldots,n-1$, if $\pi(s) < \pi(s+1)$, then $(i_{s},j_{s}) \leq^{+} (i_{s+1},j_{s+1})$, which falls into one of
two mutually exclusive cases:
\begin{align}
i_{s} \leq^{+} i_{s+1} & \mbox{ and }  j_{s}\leq^{+} j_{s+1} \mbox{ or,} \label{eqn:e1b}\\
i_{s} \leq^- i_{s+1} & \mbox{ and }  j_{s} \geq^{-} j_{s+1}.\label{eqn:e2b}
\end{align}
If $\pi(s) > \pi(s+1)$, then $(i_{s},j_{s}) \leq^{-} (i_{s+1},j_{s+1})$,
which we split into cases:
\begin{align}
i_{s} \leq^{+} i_{s+1} & \mbox{ and }  j_{s} \leq^{-} j_{s+1} \mbox{ or,} \label{eqn:e3b}\\
i_{s} \leq^- i_{s+1} & \mbox{ and }  j_{s} \geq^{+} j_{s+1}.\label{eqn:e4b}
\end{align}
We define $I_F$ to be the set of all $s$ such that either \eqref{eqn:e2b} or \eqref{eqn:e4b} holds for $F$. Notice that in both cases, $i_s \leq^{-} i_{s+1}$. Now for any $I \subset [0,n-1]$, let $S_I$ be the set of all solutions $F$ to \eqref{eq:pmb} satisfying $I_F = I$.

For any particular $I\subset [0,n-1]$, form the poset $P_{I}$ of the elements
$0,\pm 1,\pm 2,\ldots,\pm n$ by $\pi(s) <_{P_{I}} \pi(s+1)$ if $s \notin I$,
$\pi(s) >_{P_{I}} \pi(s+1)$ if $s\in I$, where we extend all our relations by the symmetry property of type B posets. We form a zig-zag poset of $n$ elements labeled consecutively by $0,\pi(1),
\pi(2),\ldots,\pi(n)$ with downward zigs corresponding to the
elements of $I$. So if $\pi = (-3,1,-2)$ and $I = \{ 0,2\}$, then our type B poset $P_I$ is \[ 2 >_{P_I} -1 <_{P_I} 3 >_{P_I} 0 >_{P_I} -3 <_{P_I} 1 >_{P_I} -2.\]

For any solution $F$ in $S_I$, let $f: [n] \to \pm[k]'$ be defined by $f(\pi(s)) = j_s$. We will show that $f$ is an enriched $P_{I}$-partition. If $\pi(s) <_{P_{I}} \pi(s+1)$ and $\pi(s) < \pi(s+1)$
in $\mathbb{Z}$, then \eqref{eqn:e1b} tells us that $f(\pi(s)) =
j_{s} \leq^{+} j_{s+1} = f(\pi(s+1))$. If $\pi(s) <_{P_{I}} \pi(s+1)$
and $\pi(s) > \pi(s+1)$ in $\mathbb{Z}$, then \eqref{eqn:e3b} tells
us that $f(\pi(s)) = j_{s} \leq^{-} j_{s+1} = f(\pi(s+1))$. If $\pi(s)
>_{P_{I}} \pi(s+1)$ and $\pi(s) < \pi(s+1)$ in $\mathbb{Z}$,
then \eqref{eqn:e2b} gives us that $f(\pi(s))= j_{s} \geq^{-} j_{s+1} =
f(\pi(s+1))$. If $\pi(s) >_{P_{I}} \pi(s+1)$ and $\pi(s) >
\pi(s+1)$ in $\mathbb{Z}$, then \eqref{eqn:e4b} gives us that
$f(\pi(s)) = j_{s} \geq^{+} j_{s+1} = f(\pi(s+1))$. In other words, we
have verified that $f$ is a $P_{I}$-partition. So for any particular solution in $S_I$, the $n$-tuple $(j_1, \ldots, j_n)$ can be thought of as an enriched $P_{I}$-partition.

Conversely, any enriched $P_{I}$-partition $f$ gives a solution in $S_I$ since if $j_s = f(\pi(s))$, then \[(( i_1,j_1),\ldots,(i_n,j_n)) \in S_I\] if and only if $0 \leq i_1 \leq \cdots \leq i_n \leq l$ and $i_s \leq^{-} i_{s+1}$ for all $s \in I$, $i_s \leq^+ i_{s+1}$ for $s \notin I$. We can therefore turn our attention to counting enriched $P_{I}$-partitions, and the remainder of the argument follows the proof of Theorem \ref{thm:ges}.
\end{proof}

\section{Quasisymmetric functions and peak algebras}\label{sec:qsym}

We will now make a connection between enriched $P$-partitions and quasisymmetric functions. Gessel \cite{Gessel} used ordinary $P$-partitions as a guide for the study of the quasisymmetric functions and Stembridge \cite{Stembridge} used enriched $P$-partitions in a similar way. We will review some of Stembridge's results, and then proceed analogously for left enriched $P$-partitions and type B enriched $P$-partitions. There are many applications for quasisymmetric functions that we will not discuss here; we primarily lay the groundwork for future studies. However, in section \ref{sec:coalgebras} we discuss the implications for subalgebras of the group algebra.

\subsection{Generating functions}

Recall that a quasisymmetric function is one for which the coefficient of $z_{i_1}^{\alpha_1} z_{i_2}^{\alpha_2} \cdots z_{i_k}^{\alpha_k}$ is the same for all fixed tuples of integers $(\alpha_1, \alpha_2, \ldots, \alpha_k)$ and all $i_1 < i_2 < \cdots < i_k$. There are two common bases for the ring of quasisymmetric functions. For any subset $S = \{ s_1 < s_2 < \cdots < s_{k-1} \}$ of $[n]$, define the monomial quasisymmetric functions, $M_S$, and the fundamental quasisymmetric functions, $F_S$:
\begin{align*}
M_S & =  \sum_{i_1 < i_2 < \cdots < i_k} z_{i_1}^{s_1} z_{i_2}^{s_2 - s_1} \cdots z_{i_k}^{n-s_{k-1}}\\
 & =  \sum_{i_1 < i_2 < \cdots < i_k} z_{i_1}^{\alpha_1} z_{i_2}^{\alpha_2} \cdots z_{i_k}^{\alpha_k}.\\
F_S & =  \sum_{S \subset T \subset [n-1]} M_T \\
 & =  \sum_{\substack{ i_1 \leq i_2 \leq \cdots \leq i_n \\ s \in S \Rightarrow i_s < i_{s+1} }} \prod_{s=1}^{n} z_{i_s}\\
 & =  \Gamma(\pi),
\end{align*}
where $\Gamma(\pi)$ is the generating function for the ordinary $P$-partitions of a permutation $\pi$ with descent set $S$. Notice that if we have the generating function for a permutation it is easy to recover its order polynomial by specializing:
\[ \Omega(\pi;m) = \Gamma(\pi)(1^m),\] where $\Gamma(\pi)(1^m)$ means that we set $z_1 = \cdots = z_m = 1$, and $z_r = 0$ for $r > m$.

The functions $M_S$ (or $F_S$), taken over all subsets $S \subset [n-1]$, span the quasisymmetric functions of degree $n$, denoted $\mathcal{Q}sym_n$. We define the ring of quasisymmetric functions by $\mathcal{Q}sym := \bigoplus_{n \geq 0} \mathcal{Q}sym_n$.

Define the generating function for enriched $P$-partitions $f: P \to \mathbb{P}'$ by \[ \dd(P) = \sum_{f \in \mathcal{E}(P)}\prod_{i=1}^{n} z_{|f(i)|}.\] Then clearly $\dd(P)$ is a quasisymmetric function. Whenever we like, we can write \[\Omega'(P;m) = \dd(P)(1^m),\] and by the fundamental Lemma \ref{lem:FLEPP}, we have that \[\dd(P) = \sum_{\pi \in \mathcal{L}(P) } \dd(\pi). \] From \cite{Stembridge} we see that the generating function for enriched $\pi$-partitions depends on the peak set, and moreover it can be written as a sum of generating functions for ordinary $P$-partitions.

\begin{thm}[Stembridge \cite{Stembridge}, Proposition 2.2]
We have the following equality: \[ \dd(\pi) = \sum_{\substack{ E \subset [n-1]
\\ \Pe(\pi) \subset E \cup (E+1) } } {\kern -10pt} 2^{|E|+1}M_{E}. \]
\end{thm}
For any subset of the integers $S$, define the set $S+1 = \{ s+1 \mid s \in S \}$. For any sets $S$ and $T$, let $S \vartriangle T = (S \cup T)\setminus (S \cap T)$ denote the symmetric difference of sets.
\begin{thm}[Stembridge \cite{Stembridge}, Proposition 3.5]
We have the following equality:
\begin{equation}\label{eq:peakdesgf}
 \dd(\pi) = 2^{\pe(\pi)+1} {\kern -10pt} \sum_{\substack{ D \subset [n-1] \\ \Pe(\pi) \subset D \vartriangle (D + 1)}} {\kern -10pt} F_D.
\end{equation}
\end{thm}

Proofs of these Theorems can be found in \cite{Stembridge}, but they also follow as special cases of our results. We point out that equation \eqref{eq:peakdesgf} specializes to \[\Omega'(\pi;m) = 2^{\pe(\pi)+1} {\kern -10pt} \sum_{\substack{ D \subset [n-1] \\ \Pe(\pi) \subset D \vartriangle (D + 1)}} {\kern -10pt} \Omega(D;m),\] as used in the proof of Theorem \ref{thm:epgf}.

For interior peak sets $S$, let $K_S$ be the quasisymmetric function defined by \[ K_{\Pe(\pi)} := \dd(\pi).\] Let $\mathbf{\Pi}_n$ denote the space of quasisymmetric functions spanned by the $K_S$, where $S$ runs over all interior peak sets of $[n-1]$. Stembridge \cite{Stembridge} defined the ``algebra of peaks" as $\mathbf{\Pi} := \bigoplus_{n\geq 0} \mathbf{\Pi}_n$, which is a graded subring of $\mathcal{Q}sym$. He proved that the functions $K_S$ are linearly independent, and so the rank of $\mathbf{\Pi}_n$ is the Fibonacci number $f_{n-1}$ (the number of distinct interior peak sets), defined by $f_0 = f_1 = 1$ and $f_n = f_{n-1} + f_{n-2}$ for $n \geq 2$.

Chow \cite{Chow} related ordinary type B $P$-partitions to type B quasisymmetric functions, and we will now discuss how the type B quasisymmetric functions relate to left enriched $P$-partitions and type B enriched $P$-partitions.  Let $S$ be any subset of $[0,n]$, $S = \{s_1 < s_2 < \cdots < s_{k-1}\}$. For fixed $n$, define the \emph{monomial} and \emph{fundamental quasisymmetric functions of type B} to be
\begin{align*}
N_S & =  \sum_{0 < i_2 < \cdots < i_k} z_{0}^{s_1} z_{i_2}^{s_2 - s_1} \cdots z_{i_k}^{n-s_{k-1}}\\
 & =  \sum_{0 < i_2 < \cdots < i_k} z_{0}^{\alpha_1} z_{i_2}^{\alpha_2} \cdots z_{i_k}^{\alpha_{k}},\\
L_S & =  \sum_{S \subset T \subset [0,n-1]} N_T\\
 & =  \sum_{\substack{ 0 \leq i_1 \leq i_2 \leq \cdots \leq i_n \\ s \in S \Rightarrow i_s < i_{s+1} }} \prod_{s=1}^{n} z_{i_s}\\
& =  \Gamma_B(\pi),
\end{align*}
where $\Gamma_B(\pi)$ is the generating function for the ordinary type B $P$-partitions of any signed permutation $\pi$ with descent set $S$. Again, we can specialize: \[\Omega_B(\pi;m) = \Gamma_B(\pi)(1^{m+1}).\] Here $\Gamma_B(\pi)(1^{m+1})$ means that we set $z_0 = z_1 = \cdots = z_m = 1$, and $z_r = 0$ for $r > m$.

The functions $N_S$ (or $L_S$), taken over subsets $S \subset [0,n-1]$, form a basis for the type B quasisymmetric functions of degree $n$, denoted $\mathcal{BQ}sym_n$. The ring of type B quasisymmetric functions is $\mathcal{BQ}sym := \bigoplus_{n \geq 0} \mathcal{BQ}sym_n$.

Define the generating function for left enriched $P$-partitions $f: P \to \mathbb{P}^{(\ell)}$, and type B enriched $P$-partitions $f: P \to \mathbb{Z}'$,
\begin{align*}
\dd^{(\ell)}(P) & =  \sum_{f \in \mathcal{E}^{(\ell)}(P)} \prod_{i =1}^{n} z_{|f(i)|},\\
\dd_B(P) & =  \sum_{f \in \mathcal{E}_B(P)} \prod_{i =1}^{n} z_{|f(i)|}.
\end{align*} We have
\begin{align*}
\lom(P;m) & =  \dd^{(\ell)}(P)(1^{m+1}), \\
\Omega_B'(P;m) & =  \dd_B(P)(1^{m+1}),
\end{align*} and the fundamental lemma gives that
\begin{align*}
\dd^{(\ell)}(P) & =  \sum_{\pi \in \mathcal{L}(P)} \dd^{(\ell)}(\pi),\\
\dd_B(P) & =  \sum_{\pi \in \mathcal{L}_B(P)} \dd_B(\pi).
\end{align*}
We can relate $\dd^{(\ell)}(\pi)$ and $\dd_B(\pi)$ to the monomial and fundamental quasisymmetric functions of type B. Notice that for a permutation $\pi \in \mathfrak{S}_n \subset \mathfrak{B}_n$, left peaks coincide with the type B peaks. Therefore we can view left enriched $P$-partitions as a special case of type B $P$-partitions. Furthermore, since enriched $P$-partitions are simply those left enriched $P$-partitions that are nonzero, we have \[ \dd(P)(z_1, z_2, \ldots) = \dd^{(\ell)}(P)(0,z_1,z_2,\ldots ),\] so the results for $\dd(P)$ can be obtained from our results for $\dd^{(\ell)}(P)$ by setting $z_0 = 0$.

\begin{thm}\label{thm:mon}
We have the following equation:
\[
\dd_B(\pi) = \sum_{ \substack{ E \subset [0,n-1] \\ \Pe_B(\pi) \subset E \cup (E+1) \\ \pi(1) < 0 \Rightarrow 0 \in E }} {\kern -10pt} 2^{|E|} N_{E}.
\]
\end{thm}

\begin{cor}
The function $\dd_B(\pi)$ depends only on the peak set of $\pi$ and the sign of $\pi(1)$.
\end{cor}

\begin{cor}
We have the following equation:
\[
\dd^{(\ell)}(\pi) = \sum_{ \substack{ E \subset [0,n-1] \\ \lPe(\pi) \subset E \cup (E+1) }} {\kern -10pt} 2^{|E|} N_{E}.
\]
\end{cor}

For any signed permutation $\pi$ with $\Pe_B(\pi) = S$, let the pair $S' = (\varsigma(\pi), S)$ be the \emph{sign-peak set} of $\pi$, where we recall $\varsigma(\pi) = 0$ if $\pi(1) > 0$, $\varsigma(\pi) = 1$ if $\pi(1) < 0$. We say $\varsigma$ is the \emph{sign} of $S'$. If $S$ is any valid type B peak set, then in general we have two associated sign-peak sets, $(0,S)$ and $(1,S)$. But if position 1 is a peak of $\pi$, it must be that $\pi(1)>0$. In other words, if $1 \in S$, the pair $(0,S)$ is the only valid sign-peak set. For fixed $n$, we can see that the number of sign-peak sets is the Fibonacci number $f_{n+1}$.

We define the functions $K_{S'}$ by \[K_{(\varsigma(\pi), \Pe_B(\pi))} := \dd_B(\pi).\]
Note that for a permutation $\pi \in \mathfrak{S}_n$, $\dd^{(\ell)}(\pi) = K_{(0,\lPe(\pi))}$.
\begin{thm}\label{thm:fun}
We have the following equation:
\[ \dd_B(\pi) = 2^{\pe_B(\pi) + \varsigma(\pi)} {\kern -10pt} \sum_{ \substack{ D \subset [0,n-1] \\ \Pe_B(\pi) \subset D \vartriangle (D+1) \\ \pi(1) < 0 \Rightarrow 0 \in D }} {\kern -10pt} L_D\]
In other words,
\begin{align*}
K_{(0,S)} & =  2^{|S|} {\kern -10pt} \sum_{\substack{ D \subset [0,n-1] \\ S \subset D \vartriangle (D+1) } } {\kern -10pt} L_D \\
K_{(1,S)} & =  2^{|S|+1} {\kern -10pt} \sum_{\substack{ 0\in D \subset [0,n-1]  \\ S \subset D \vartriangle (D+1) } } {\kern -10pt} L_D
\end{align*}
\end{thm}
This theorem specializes to the formulas
\begin{align*}
\lom(\pi;m) & =  2^{\lpe(\pi)} {\kern -10pt} \sum_{\substack{ D \subset [0,n-1] \\ \lPe(\pi) \subset D \vartriangle (D+1) } } {\kern -10pt} \Omega_B(D;m), \\
\Omega'_B(\pi;m) &=  2^{\pe_B(\pi) + \varsigma(\pi)} {\kern -10pt} \sum_{ \substack{ D \subset [0,n-1] \\ \Pe_B(\pi) \subset D \vartriangle (D+1) \\ \pi(1) < 0  \Rightarrow 0 \in D }} {\kern -10pt} \Omega_B(D;m),
\end{align*}
as used in the proofs of Theorems \ref{thm:eepgf} and \ref{thm:bepgf}.

Let $\mathbf{\Pi}_{B,n}$ denote the span of the $K_{S'}$, where $S'$ ranges over all sign-peak sets of $[n-1]$. It is not hard to see that the $K_{S'}$ are linearly independent, and so $\mathbf{\Pi}_{B,n}$ has rank $f_{n+1}$. If we define the type B algebra of peaks, $\mathbf{\Pi}_B := \bigoplus_{n \geq 0} \mathbf{\Pi}_{B,n}$, then we can see it is a subring of $\mathcal{BQ}sym$, as an argument identical to that of \cite{Stembridge} Theorem 3.1 shows.

Similarly, let $\mathbf{\Pi}^{(\ell)}_n$ denote the span of all $K_{S'}$, where $S'$ ranges over those sign-peak sets with sign zero. Then $\mathbf{\Pi}^{(\ell)}_n$ has rank $f_n$ and the left algebra of peaks, $\mathbf{\Pi}^{(\ell)}:= \bigoplus_{n\geq 0} \mathbf{\Pi}^{(\ell)}_n$, is a subring of $\mathbf{\Pi}_B$.

Now we prove the above theorems.

\begin{proof}[Proof of Theorem \ref{thm:mon}]
Clearly, we can expand $\dd_B(\pi)$ as a sum of $N_E$ with nonnegative coefficients. Fix $E = \{ a_1 < a_2 < \cdots < a_k \} \subset [0,n-1]$. The coefficient of $N_E$ is equal to the coefficient of $z_{0}^{\alpha_1} z_{1}^{\alpha_2} \cdots z_{k}^{\alpha_{k+1}}$, where $\alpha_i = a_i - a_{i-1}$, with $a_0 = 0$, $a_{k+1} = n$. This coefficient is equal to the number of type B enriched $\pi$-partitions $f$ such that
\begin{equation}\label{eq:tuple}
(|f(\pi(1))|, |f(\pi(2))|, \ldots, |f(\pi(n))|) = ( \underbrace{0,0,\ldots, 0}_{\alpha_1}, \underbrace{ 1,1,\ldots , 1}_{\alpha_2}, \ldots, \underbrace{k,k,\ldots, k}_{\alpha_{k+1}} ).
\end{equation}
In other words, all zeros until position $a_1$, all ones from position $a_1 + 1$ to $a_2$, and so on. Notice that within each block of numbers, the signed permutation $\pi$ must satisfy one of three conditions. It must either be always increasing, always decreasing, or decreasing then increasing. It cannot be increasing then decreasing, since then it would have a peak. Say there is a peak in position $i$. Then $f(\pi(i-1)) \leq^{+} f(\pi(i)) \leq^{-} f(\pi(i+1))$, so $|f(\pi(i-1))| < |f(\pi(i+1))|$. Therefore the only possible positions for peaks are $a_1$ or $a_1 + 1$, $a_2$ or $a_2 + 1, \ldots, a_k$ or $a_k +1$. In other words, the coefficient of $N_E$ is nonzero only if $\Pe_B(\pi) \subset E \cup (E+1)$.

Now given that $\Pe_B(\pi) \subset E \cup (E+1)$, we will determine the coefficient of $N_E$. Within each nonzero block of numbers, we have ($i \neq 0$): \[ ( |f(\pi(a_{i} + 1))|, |f(\pi(a_i + 2))|, \ldots, |f(\pi(a_{i+1}))| ) = ( i,i, \ldots, i).\] We claim there are exactly two possibilities for $f$ in every such block. If $\pi$ is increasing over this interval, then \[f(\pi( a_i + 1)) \leq^{+} f(\pi( a_i + 2)) \leq^{+} \cdots \leq^{+} f(\pi(a_{i+1})),\] so then $f(\pi(a_i + 1)) = \pm i$, and all others equal $+i$. If $\pi$ is decreasing over the entire interval, then \[f(\pi( a_i + 1)) \leq^{-} f(\pi( a_i + 2)) \leq^{-} \cdots \leq^{-} f(\pi(a_{i+1})),\] so then $f(\pi(a_{i + 1})) = \pm i$, and all others equal $-i$. The third case has $\pi$ decreasing, then increasing. Suppose $\pi(j-1) > \pi(j) < \pi(j+1)$ with $a_{i}+1 < j < a_{i+1}$. Then \[f(\pi( a_i + 1)) \leq^{-} \cdots \leq^{-} f(\pi(j)) \leq^{+} \cdots \leq^{+} f(\pi(a_{i+1})),\] so then $f(\pi(j)) = \pm i$, everything in the block to its left is equal to $-i$, while everything to its right is $+i$.

In total, there are $k$ blocks where some choice can be made, and so there are $2^k = 2^{|E|}$ such $f$. We are almost finished with the proof. The final observation to make is that if $\pi(1) < 0$, then $0 \leq^{-} f(\pi(1))$, which means that $|f(\pi(1))| > 0$, and so there can be no leading zeros in the $n$-tuple \eqref{eq:tuple}. In other words, it must be that $\alpha_1 = a_1 = 0$, and the theorem is proved.
\end{proof}

\begin{proof}[Proof of Theorem \ref{thm:fun}]
Let us suppose $\pi(1) > 0$ and expand the following in terms of the $N_E$:
\begin{equation}\label{eq:expand1}
\sum_{ \substack{ D \subset [0,n-1] \\ \Pe_B(\pi) \subset D \vartriangle (D+1)} } {\kern -10pt} L_D = \sum_{\substack{ D \subset [0,n-1] \\ \Pe_B(\pi) \subset D \vartriangle (D+1) } } \sum_{D \subset E \subset [0,n-1] } N_E
\end{equation}
The coefficient of $N_E$ in \eqref{eq:expand1} is $\# \{ D\subset E : \Pe_B(\pi) \subset D \vartriangle (D+1) \}$, which is clearly zero unless $\Pe_B(\pi) \subset  E \cup (E+1)$. Now for any such $E$, we have several cases. We have $j, j+1$ both in $E$ and $j+1$ in $\Pe_B(\pi)$, exactly one of $j$ or $j+1$ in $E$ and $j+1$ in $\Pe_B(\pi)$, or $j$ in $E$ but neither $j$ nor $j+1$ is in $\Pe_B(\pi)$.

If both $j$ and $j+1$ are in $E$, $j+1$ in $\Pe_B(\pi)$, then $D$ can have exactly one of $j$ or $j+1$. If exactly one of $j$ or $j+1$ in $E$ and $j+1$ in $\Pe_B(\pi)$, then $D$ must have whichever $E$ has. If $j$ is in $E$ but neither $j$ nor $j+1$ is in $\Pe_B(\pi)$, then $j$ is free to be in $D$ or not. For example, if $E = \{ 1,2, 3, 5,6,8,9\}$ and $\Pe_B(\pi) = \{ 2,4,5,7\}$, then $D$ must have $3,5,6$, and exactly one of $1$ or $2$. It is free to contain $8$ or $9$ (or not). There are $8 = 2^3$ such $D$, and it should be clear from this example that there should be $2^{|E| - \pe_B(\pi)}$ choices in general. Therefore we have the coefficient of $N_E$ in \eqref{eq:expand1} is $2^{|E|-\pe_B(\pi)}$, and
\begin{align*}
2^{\pe_B(\pi)} {\kern -10pt} \sum_{\substack{ D \subset [0,n-1] \\ \Pe_B(\pi) \subset D \vartriangle (D+1) } } {\kern -10pt} L_D & =  2^{\pe_B(\pi)} {\kern -5pt} \sum_{\Pe_B(\pi) \subset E \cup (E+1)} {\kern -5pt} 2^{|E|-\pe_B(\pi)} N_E \\
 & =  \sum_{ \Pe_B(\pi) \subset E \cup (E+1) } {\kern -5pt} 2^{|E|} N_E \\
 & =  \dd_B(\pi) \qquad \mbox{ (by Theorem \ref{thm:mon}). }
\end{align*}

Now let us suppose $\pi(1) < 0$ and expand the following in terms of the $N_E$:
\begin{equation}\label{eq:expand2}
\sum_{ \substack{ 0 \in D \subset [0,n-1] \\ \Pe_B(\pi) \subset D \vartriangle (D+1) } } {\kern -10pt} L_D = \sum_{\substack{ 0\in D \subset [0,n-1] \\ \Pe_B(\pi) \subset D \vartriangle (D+1) } } \sum_{0 \in D \subset E \subset [0,n-1] } N_E.
\end{equation}

The argument is identical to the $\pi(1) > 0$ case, except that now we have both $E$ and $D$ must contain 0, so we have one fewer choice to make in selecting possible subsets $D$. Specifically, the coefficient of $N_E$ in \eqref{eq:expand2} is $2^{|E|-\pe_B(\pi)-1}$. Therefore we have,

\begin{align*}
2^{\pe_B(\pi)+1} {\kern -10pt} \sum_{\substack{ 0\in D \subset [0,n-1] \\ \Pe_B(\pi) \subset D \vartriangle (D+1) } } {\kern -10pt} L_D & =  2^{\pe_B(\pi)+1} {\kern -10pt} \sum_{\substack{ \Pe_B(\pi) \subset E \cup (E+1) \\ 0 \in E }} {\kern -5pt} 2^{|E|-\pe_B(\pi)-1} N_E \\
 & =  \sum_{\substack{ \Pe_B(\pi) \subset E \cup (E+1) \\ 0 \in E }} {\kern -5pt} 2^{|E|} N_E \\
 & =  \dd_B(\pi) \qquad \mbox{ (by Theorem \ref{thm:mon}). }
\end{align*}
The proof is complete.
\end{proof}

\subsection{Coalgebra structures}\label{sec:coalgebras}

Let $X = \{ x_1, x_2, \ldots \}$ and $Y = \{ y_1, y_2, \ldots \}$ be two sets of commuting indeterminates. Define the set $XY = \{ (x,y) : x \in X, y \in Y\}$. Then we define the bipartite generating function, \[ \Gamma(\pi)(XY) = \sum_{\substack{ (i_1,j_1) \leq (i_2, j_2) \leq \cdots \leq(i_n,j_n) \\ \pi(s) > \pi(s+1) \Rightarrow (i_s,j_s) < (i_{s+1}, j_{s+1}) }} {\kern -10pt} x_{i_1}\cdots x_{i_n} y_{i_1} \cdots y_{i_n}. \] In \cite{Gessel}, Gessel proved the following theorem (which we specialized for our Theorem \ref{thm:ges}).

\begin{thm}[Gessel \cite{Gessel}]\label{thm:gesqsym}
We have the following equation:
\begin{equation}\label{eq:internal}
\Gamma(\pi)(XY) = \sum_{\sigma\tau = \pi} \Gamma(\sigma)(X) \Gamma(\tau)(Y).
\end{equation}
\end{thm}

One implication of equation \eqref{eq:internal}, established in \cite{Gessel}, is that the space of quasisymmetric functions of degree $n$ forms a coalgebra dual to Solomon's descent algebra. Specifically, if $\pi$ is any permutation with $\Des(\pi) = K$, let $a_{I,J}^{K}$ denote the number of pairs of permutations $(\sigma, \tau) \in \mathfrak{S}_n \times \mathfrak{S}_n$ with $\Des(\sigma) = I$, $\Des(\tau) = J$, and $\sigma\tau = \pi$. Then \eqref{eq:internal} may be restated as the coproduct $\mathcal{Q}sym_n \to \mathcal{Q}sym_n \otimes \mathcal{Q}sym_n$: \[ F_K \mapsto \sum_{I,J \subset [n-1]} a_{I,J}^{K} F_I \otimes F_J.\] The Solomon descent algebra is then isomorphic to the dual space $\mathcal{Q}sym_n^{*}$ with multiplication \[ u_I \ast u_J = \sum_K a_{I,J}^K u_K,\] where we recall from the introduction that $u_I$ is the sum of all permutations with descent set $I$.

The proof of Theorem \ref{thm:gesqsym} is the same as the proof given for Theorem \ref{thm:ges}, except that we take the lexicographic order on all of $\mathbb{P}\times \mathbb{P}$, rather than on a finite grid like $[l]\times[k]$. Similarly, we can extend the proofs for Theorems \ref{thm:interior}, \ref{thm:left}, \ref{thm:peakideal}, \ref{thm:interiordescent}, and \ref{thm:peakalg2}. With some of these theorems we are treading in the territory of type B quasisymmetric functions, so we define $X_0 = X \cup \{x_0\}$ for any set $X$.

\begin{thm}\label{thm:interiorgf}
We have the following equation:
\begin{equation}\label{eq:interiorgf}
\dd(\pi)(XY) = \sum_{\sigma\tau = \pi} \dd(\sigma)(X) \dd(\tau)(Y).
\end{equation}
\end{thm}

Here we extend the proof of Theorem \ref{thm:interior} by considering the up-down order on $\mathbb{P}'\times\mathbb{P}'$. Let $\mathfrak{P}_n$ denote the span of $v_I$, sums of permutations with the same set of interior peaks. Then \eqref{eq:interiorgf} tells us that $\mathbf{\Pi}_n$ is the coalgebra dual to $\mathfrak{P}_n$, with comultiplication \[ K_U \mapsto \sum_{S,T} c_{S,T}^U K_S \otimes K_T,\] where the sum ranges over all pairs of interior peak sets and for any permutation $\pi$ with $\Pe(\pi) = U$, $c_{S,T}^U$ is the number of pairs of permutations $(\sigma, \tau)$ with $\Pe(\sigma) = S$, $\Pe(\tau) = T$, and $\sigma\tau = \pi$.

\begin{thm}\label{thm:leftgf}
We have the following equation:
\begin{equation}\label{eq:leftgf}
\dd^{(\ell)}(\pi)(X_0Y_0) = \sum_{\sigma\tau = \pi} \dd^{(\ell)}(\sigma)(X_0) \dd^{(\ell)}(\tau)(Y_0).
\end{equation}
\end{thm}

Here we modify the proof of Theorem \ref{thm:left} to take the up-down order on $\mathbb{P}^{(\ell)}\times\mathbb{P}^{(\ell)}$. We let $\mathfrak{P}^{(\ell)}_n$ denote the span of \[v^{(\ell)}_I := \sum_{\lPe(\pi) = I} \pi,\] sums of permutations with the same set of left peaks. Then \eqref{eq:leftgf} implies that $\mathbf{\Pi}^{(\ell)}_n$ is the coalgebra dual to $\mathfrak{P}^{(\ell)}_n$, with comultiplication \[ K_{(0,U)} \mapsto \sum_{S,T} d_{S,T}^{U} K_{(0,S)} \otimes K_{(0,T)},\] where the sum ranges over all pairs of left peak sets (sign-peak sets with sign zero) and for any $\pi$ such that $\lPe(\pi) = U$, $d_{S,T}^{U}$ is the number of pairs of permutations $(\sigma, \tau)$ with $\lPe(\sigma) = S$, $\lPe(\tau) = T$, and $\sigma\tau = \pi$.

We remark that while much is already known about $\mathfrak{P}_n$ and $\mathfrak{P}^{(\ell)}_n$ from \cite{AguiarBergeronNyman}, the structure constants $c_{S,T}^{U}$, $d_{S,T}^{U}$, lacked the combinatorial description we provide here.

\begin{thm}\label{thm:peakalg2gf}
We have the following equation:
\begin{equation}\label{eq:peakalg2gf}
\dd_B(\pi)(X_0Y_0) = \sum_{\sigma\tau = \pi} \dd_B(\sigma)(X_0) \dd_B(\tau)(Y_0).
\end{equation}
\end{thm}

Here we take the up-down order on $\mathbb{Z}'\times\mathbb{Z}'$, and the proof goes through as for Theorem \ref{thm:peakalg2}. Let
\[
v_{I'}  := \sum_{(\varsigma(\pi), \Pe_B(\pi)) = I'} \pi,
\]
and let $\mathfrak{P}_{B,n}$ denote the span of the $v_{I'}$. Then Theorem \ref{thm:peakalg2gf} shows that $\mathbf{\Pi}_{B,n}$ is the coalgebra dual to $\mathfrak{P}_{B,n}$. Indeed, we can interpret \eqref{eq:peakalg2gf} as giving the map
\[
K_{U'} \mapsto \sum_{S',T'} b_{S',T'}^{U'} K_{S'} K_{T'},
\]
where the sum ranges over all pairs of sign-peak sets and for any signed permutation $\pi$ with $(\varsigma(\pi), \Pe_B(\pi)) = U'$, $b_{S',T'}^{U'}$ is the number of pairs of signed permutations $(\sigma, \tau)$ such that $(\varsigma(\sigma), \Pe_B(\sigma)) = S'$, $(\varsigma(\tau), \Pe_B(\tau)) = T'$, and $\sigma\tau = \pi$. The multiplication in $\mathfrak{P}_{B,n}$ is given explicitly by \[ v_{S'} \ast v_{T'} = \sum_{U'} b_{S',T'}^{U'} v_{U'},\] the sum taken over all sign-peak sets $U'$.

The generalizations of Theorems \ref{thm:peakideal} and \ref{thm:interiordescent} tell us that $\mathfrak{P}_n$ is a left ideal in $\Sol(A_{n-1})$, and a two-sided ideal in $\mathfrak{P}^{(\ell)}_n$.

\begin{thm}\label{thm:peakidealgf}
We have the following equations:
\begin{align*}
\dd(\pi)(XY_0) & = \sum_{\sigma\tau = \pi} \dd^{(\ell)}(\sigma)(Y_0) \dd(\tau)(X),\\
\dd(\pi)(X_0Y) & = \sum_{\sigma\tau = \pi} \dd(\sigma)(Y) \dd^{(\ell)}(\tau)(X_0).
\end{align*}
\end{thm}

Here we take the up-down order on $\mathbb{P}'\times\mathbb{P}^{(\ell)}$ or $\mathbb{P}^{(\ell)}\times\mathbb{P}'$.

\begin{thm}\label{thm:interiordescentgf}
We have the following equation:
\[
\dd(\pi)(XY) = \sum_{\sigma\tau = \pi} \dd(\sigma)(X) \Gamma(\tau)(Y).
\]
\end{thm}

For this proof we take the lexicographic order on $\mathbb{P}\times\mathbb{P}'$, with the convention that $(i,j) \leq^+ (i',j')$ if $i < i'$ or $i = i'$ and $j \leq^+ j'$ and similarly for $\leq^-$.

We finish with a list of some subalgebras that \emph{do not} exist in general. For the symmetric group, the sums of permutations with the same set of right peaks do not form an algebra, nor do the sums of permutations with the same set of exterior peaks. For the hyperoctahedral group, the sums of permutations with the same sign on $\pi(1)$ and the same exterior peak set does not work, and neither does collecting signed permutations with the same number of exterior peaks and sign on $\pi(1)$.

\end{document}